\documentclass[a4paper,11pt]{article}

\renewcommand{\baselinestretch}{2}
\usepackage{fontawesome5} % <============================================
\usepackage{graphicx}

\usepackage{lscape}

\usepackage{diagbox}

\usepackage[T1]{fontenc}
\usepackage[latin1]{inputenc}
\usepackage{latexsym} %Simbolos adicionales de LaTeX
\usepackage{amssymb}
\usepackage{fancyvrb}
\usepackage{algorithmic}

\usepackage[T1]{fontenc}
\usepackage[latin1]{inputenc}
\usepackage[english]{babel}
\usepackage{latexsym} %Simbolos adicionales de LaTeX
\usepackage{amssymb}
\usepackage{amsmath}
\usepackage{array}
\usepackage{bbding} % paquete para estrellas de simbolos
\usepackage{mathpazo} %Font type: Palatino

\usepackage{amsthm}
\usepackage{xcolor}

\usepackage{apacite}
\usepackage{natbib} %For BibTeX

\usepackage{setspace}
\usepackage[font=small]{caption}
\usepackage{rotating}

\usepackage{subcaption}
\usepackage{tikz}
\usetikzlibrary{backgrounds}

\usetikzlibrary{shapes.geometric}
\usetikzlibrary{arrows,automata}

\usepackage{multirow}

\usepackage{color}
\usepackage{colortbl}

\usepackage{pgf}
\usepackage{mathrsfs}
\usetikzlibrary{arrows}

\usepackage{algorithmic}
\usepackage{longtable}
\usepackage{arydshln}
\usepackage{hyperref}
\usepackage{pdflscape}
\allowdisplaybreaks
\usepackage[margin=1 in]{geometry}

\theoremstyle{plain}

\newtheorem{theorem}{Theorem}[section]
\newtheorem{proposition}[theorem]{Proposition}

\newtheorem{definition}[theorem]{Definition}
\newtheorem{example}[theorem]{Example}

\theoremstyle{definition}

\theoremstyle{remark}

\newtheorem{remark}{Remark}[section]

\newcommand{\abs}[1]{\left|#1\right|}

\usepackage{floatrow}
\newfloatcommand{capbtabbox}{table}[][\FBwidth]

\usepackage{blindtext}

\newcommand{\addQEDstyle}[2]{\AtBeginEnvironment{#1}{\pushQED{\qed}\renewcommand{\qedsymbol}{#2}}\AtEndEnvironment{#1}{\popQED}}
\addQEDstyle{example}{$\triangle$}

\usepackage{siunitx}

\newcolumntype{R}[2]{%
    >{\adjustbox{angle=#1,lap=\width-(#2)}\bgroup}%
    l%
    <{\egroup}%
}
\allowdisplaybreaks 

\begin{document}

\title{\textbf{On the influence of dependent features in classification problems: a game-theoretic perspective}}

\author{Laura Davila-Pena$^{1,2}$, Alejandro Saavedra-Nieves$^3$, Balbina Casas-M\'endez$^4$}

\date{\footnotesize \emph{
		$^1$\underline{Corresponding author}. Department of Analytics, Operations and Systems, Kent Business School, University of Kent, CT2 7PE Canterbury, UK. ORCID: \href{https://orcid.org/0000-0003-2175-2546}{0000-0003-2175-2546}. \texttt{l.davila-pena@kent.ac.uk}\\
  $^2$MODESTYA Research Group,  Department of Statistics, Mathematical Analysis and Optimization, Faculty of Mathematics, University of Santiago de Compostela, Campus Vida, 15782 Santiago de Compostela, Spain. \texttt{lauradavila.pena@usc.es}\\
	$^3$CITMAga, MODESTYA Research Group, Department of Statistics, Mathematical Analysis and Optimization, Faculty of Mathematics, University of Santiago de Compostela, Campus Vida, 15705 Santiago de Compostela, Spain. ORCID: \href{https://orcid.org/0000-0003-1251-6525}{0000-0003-1251-6525}. \texttt{alejandro.saavedra.nieves@usc.es} \\
        $^4$CITMAga, MODESTYA Research Group, Department of Statistics, Mathematical Analysis and Optimization, Faculty of Mathematics, University of Santiago de Compostela, Campus Vida, 15782 Santiago de Compostela, Spain. ORCID: \href{https://orcid.org/0000-0002-2826-218X}{0000-0002-2826-218X}. \texttt{balbina.casas.mendez@usc.es} 
        }}

\maketitle
%\thispagestyle{empty}

%\clearpage
%\setcounter{page}{1}
\onehalfspacing
%\doublespacing
%\thispagestyle{empty}
%	\date{\today}
\renewcommand{\baselinestretch}{1.5}

\abstract{This paper deals with a new measure of the influence of each feature on the response variable in classification problems, accounting for potential dependencies among certain feature subsets. Within this framework, we consider a sample of individuals characterized by specific features, each feature encompassing a finite range of values, and classified based on a binary response variable. This measure turns out to be an influence measure explored in existing literature and related to cooperative game theory. We provide an axiomatic characterization of our proposed influence measure by tailoring properties from the cooperative game theory to our specific context. Furthermore, we demonstrate that our influence measure becomes a general characterization of the well-known Banzhaf-Owen value for games with a priori unions, from the perspective of classification problems. The definitions and results presented herein are illustrated through numerical examples and various applications, offering practical insights into our methodologies.}

\vspace{0.5cm}


%{\bf JEL classification:} C38, C71.

\section{Introduction}\label{sec:intro}

Understanding how features influence a binary response variable becomes crucial in classification problems, with applications ranging from medical diagnosis to customer behavior analysis. In this paper, we start from a sample of individuals for which a set of characteristics that take a finite number of values has been measured. These individuals are categorized using a classifier according to a response variable that can also take a finite set of values. Our objective is to provide a model-agnostic measure of the influence of each feature on predicting the response variable in cases where certain subsets of features are dependent. 
The motivation behind this research, along with its theoretical and practical significance, stems from the absence of a game theoretic-based influence measure in the literature that accounts for feature dependencies and is both experimentally validated and theoretically grounded. We fill this gap by proposing such a measure and providing an axiomatic characterization following the Banzhaf value \citep{Banzhaf1965} approach, a game theoretical solution concept closely related to the well-known Shapley value \citep{Shapley1953}. 
The main difference between the two lies in the theoretical properties they fulfill, rendering one more suitable than the other depending on the nature of the problem under consideration \citep{feltkamp1995}. The Shapley value has recently found widespread application in the context of the interpretability of machine learning models \citep{Lundberg2017}, and both the Banzhaf and Shapley values rely on players' contributions to certain coalitions. On the other hand, the Banzhaf value has been successfully employed in various contexts to determine rankings. Recent applications include its integration into a control framework, where it offers closed-form expressions as a function of the state and establishes steady-state conditions \citep{Muros2017}. Additionally, it has been used to define a measure of political power in a Spanish autonomous community \citep{Mikel2020} and to analyze the risk of a terrorist attack by ranking terrorists within a network \citep{Algaba2024}.

In our case, the proposed measure of feature influence emerges as a generalization of the one studied in \cite{Datta2015} for scenarios with independent features, and is identified with the Banzhaf-Owen value \citep{Owen1981} for classification problems where both dependent features and the response variable are binary. %, i.e., they take only two values, quantitative or qualitative, distinct. 
Although the problem addressed is not inherently novel, this is the first time, to the best of our knowledge, that an influence measure grounded in the Banzhaf-Owen value is proposed, constituting the primary contribution of our work. The Banzhaf-Owen value is a cooperative game theory solution that extends upon the Banzhaf value and is particularly relevant in situations where players exhibit affiliations due to political, economic, or other considerations. In addition, we offer an axiomatic characterization of the influence measure for cases where the response variable is binary, implement it using the software \texttt{R}, and validate its efficacy through a variety of examples and real-world applications.

The organization of the paper is as follows.  Section~\ref{sec:lit_rev} presents a review of related literature and Section~\ref{sec:preli} offers the basic preliminaries concerning both influence measures and cooperative games. In Section~\ref{sec:main}, we provide the axiomatic characterization of our measure and Section~\ref{sec:b-o} relates our proposal with the Banzhaf-Owen value. Our methodology is validated and applied in Section~\ref{sec:results} through various numerical experiments. The paper finishes with some concluding remarks in Section~\ref{sec:conclusions}.

\section{Related literature}\label{sec:lit_rev}
%In this section, we review the most relevant related literature. 
A critical challenge for the future in modern data science and artificial intelligence is achieving adequate levels of explanation for the knowledge derived from diverse methods and techniques, ensuring their appropriate utilization \citep{Carrizosa2021}. Since 2018, the European Union has recognized the so-called ``right to explanation'' concerning algorithms designed for decision-making processes \citep{European2020}. 
\cite{Burkart2021} discuss the reasons behind the necessity for explanations of machine learning models and the various domains where such explanations are in demand, including healthcare, the automobile industry, and recommender systems, among others. 
A considerable body of literature is dedicated to the pursuit and enhancement of explanation in machine learning. A prevalent objective is to quantify the influence of the characteristics on predictions, often facilitating the selection of relevant features and enabling the analysis of model behavior with a reduced set of predictor variables. Different methodologies have been explored to address the problem of explaining machine learning models. 
One of them involves considering permutations of the sample values 
\citep{Altmann2010}. Additionally, some approaches are tailored to specific methods. For instance, within the realm of classification, \cite{Ghaddar2018} propose an iterative procedure for feature selection when making use of support vector machines (SVMs). Another suggested approach consists in using concepts from cooperative game theory, assessing each feature's contribution to prediction within any coalition of features \citep{Carrizosa2021}. 
% Version anterior, comprobar si sigo diciendo lo mismo: Another common approach makes use of cooperative game theory concepts and the idea of contribution to the prediction of each feature to any set or coalition of features (Carrizosa et al., 2021)
% Lo que dice Carrizosa  (no tengo muy claro qué quiere decir):
%  Examples of these measures[...] are based [...], inspired by cooperative games, on calculating the contribution toward the accuracy of the feature to any coalition of features
Table~\ref{tab:Literature_gen} schematically shows the references cited so far, underscoring the notable interest in the literature regarding explanations for classification models and, more broadly, for machine learning.

\begin{table}[h]
\resizebox{0.88\linewidth}{!}{
    \begin{tabular}{p{5cm}ll}
    \hline
    Context & Reference & Key point\\%[0.5em]
    \hline
    Measure of features  & \cite{Altmann2010}& Permutations of the \\
    importance &  & sample values\\[0.5em]
    Feature selection & \cite{Ghaddar2018} & Classification with SVMs\\[0.5em]
    Artificial intelligence & \cite{European2020}&Right of explanation\\[0.5em]
    Future challenges &\cite{Carrizosa2021} & Cooperative games \\[0.5em]
    Revision & \cite{Burkart2021} & Reasons/domains\\%[0.5em]
    \hline
    \end{tabular}
}
\caption{Summary of recent references on explaining classification models.}
\label{tab:Literature_gen}
\end{table}

It is also noteworthy to mention the work of \cite{Cohen2007}, who present a new feature selection algorithm on the basis of Shapley's contribution values of the features to classification precision. 
This method iteratively estimates the utility of features, facilitating subsequent forward selection or backward elimination procedures. Meanwhile, 
\cite{Kononenko2010} provide a local measure of the influence of the features in classification problems, i.e., operating at the individual instance level, and make use of the Shapley value. This procedure is extended in \cite{Strumbelj2011} to explain regression models and their predictions for individual instances. A paper close to the preceding two studies is \cite{Datta2015}. In this work, the authors explore the global influence of different features in a classification problem. Theoretically, they restrict their influence measure to scenarios where the response variable takes only two values, and it is for these specific cases that they discuss its relation with the Banzhaf value of so-called simple cooperative games. Another relevant reference within the context of interpretability in complex machine learning models leveraging game theory is due to \cite{Lundberg2017}. They introduce SHAP (SHapley Additive exPlanations), a unified framework for interpreting predictions. This framework is supported by theoretical foundations across several prediction models, with a focus on determining the importance of the features for each specific prediction. In \cite{Casalicchio2019}, alongside introducing a local measure of feature importance for individual observations and two visualization tools, the authors describe a procedure based on the Shapley value. This approach equally distributes the overall model performance among features according to their marginal contributions, enabling comparison of the importance of features across different models. \cite{Smith2021} employ the methodology proposed in \cite{Lundberg2017} to analyze a COVID-19 database collected in the early stages of the pandemic in Wuhan, China. Meanwhile, \cite{Jothi2021} concentrate on implementing a novel feature selection and data mining classifier system. According to their method, the Shapley value is set to serve as the feature selector for the data mining classifier applied to mental health data. 
\cite{Davila2022} expand upon the solution proposed in \cite{Kononenko2010} by providing a global measure of the influence of features in a classification problem, along with an axiomatic characterization. The authors apply this methodology to a sample of COVID-19 patients obtained in Spain during the first wave of the pandemic, the first four-month period of 2020. Recently, \cite{Nahiduzzaman2024} introduce a novel approach for accurately classifying three types of lung cancer together with normal lung tissue by using CT (Computed Tomography) images. The integration of SHAP into their framework enhances explanatory capabilities, offering valuable insights for decision-making and bolstering confidence in real-world lung cancer diagnoses.
Table~\ref{tab:Literature} summarizes the references where methods derived from cooperative game theory are proposed and applied to improve the interpretation of machine learning models. 

\begin{table}[ht!]
\resizebox{0.88\linewidth}{!}{
    \begin{tabular}{p{5cm}ll}
     \hline
    Context & Reference & Game theory value\\%[0.5em]
     \hline
    Selection of features & \cite{Cohen2007} & Shapley value\\[0.5em]
    Explanation of individual classifications & \cite{Kononenko2010} & Shapley value\\[0.5em]
    Explanation of regression predictions & \cite{Strumbelj2011} &    Shapley value      \\[0.5em]
    Explanation of classification of a database & \cite{Datta2015} & Banzhaf value\\[0.5em]
    Interpretation of predictions & \cite{Lundberg2017} &   Shapley value       \\[0.5em]
    Importance of features across models      & \cite{Casalicchio2019} &  Shapley value        \\[0.5em]
    Identifying mortality factors & \cite{Smith2021}&Shapley value\\[0.5em]
    Mental health data classification  & \cite{Jothi2021} &  Shapley value        \\[0.5em]
    Global explanation in classification & \cite{Davila2022} & Shapley value \\[0.5em]
    Lung cancer classification & \cite{Nahiduzzaman2024} & Shapley value\\
    %\textcolor{blue}{Explanation of regression predictions within the realm of wind power} & \textcolor{blue}{\cite{Cakiroglu2024}} & \textcolor{blue}{Shapley value}\\    %\noalign{\smallskip}\hline\noalign{\smallskip}
     \hline
    \end{tabular}
}
\caption{Summary of recent references on game theoretic methods for machine learning models.}
\label{tab:Literature}
\end{table}

According to \cite{Kononenko2010}, the main shortcoming of existing general explanation methods is that they do not consider all possible dependencies and interactions among feature values. Since then, researchers have worked to incorporate such dependencies in their proposals, although it should be noted that dependency has not been explicitly addressed from a game theory perspective.
In linear regression problems, \cite{Owen2017} analyze the global explanation of the model using Shapley's contributions of the features and address the problem of dependent features using an ANOVA decomposition approach. \cite{Giudici2021} extend this framework using the so-called Lorenz Zonoid decomposition.
\cite{Aas2021} address the problem of explaining individual predictions by treating dependent features using an extension of the Kernel SHAP method which is a computationally efficient approximation of Shapley values for the case of large problems. The method is illustrated with examples of linear and non-linear models. 

\begin{table}[h!]
\resizebox{0.88\linewidth}{!}{
    \begin{tabular}{p{5cm}ll}
     \hline
    Context & Reference & Solution concept/approach\\%[0.5em]
     \hline
    Linear regression model & \cite{Owen2017} & Shapley value/ANOVA \\[0.5em]
    Linear regression model  & \cite{Giudici2021}  & Shapley value/Lorenz Zonoids   \\[0.5em]
    Linear and non linear models & \cite{Aas2021} & Extended Kernel SHAP \\[0.5em]
    To estimate the conditional expectations                             &  \cite{Olsen2022} & Local Shapley value\\[0.5em]
    General overview                                  & \cite{Li2024} & Shapley/Owen values\\[0.5em]
    Classification  & This paper  &  Banzhaf-Owen value    \\%[0.5em]
     \hline
    \end{tabular}
}
\caption{Summary of recent references on game theoretic methods for machine learning models with dependence among features.}
\label{tab:Literature_dep}
\end{table}

In contrast to papers that use methods from the statistical domain to model dependencies between features, \cite{Olsen2022} use machine learning methods to compute local Shapley values in a regression environment. \cite{Li2024} provide a recent overview of the Shapley value as one of the main approaches from artificial intelligence to explain machine learning models, including the consideration of possible complex dependencies between features. The work reveals the potential of the Owen value (\citealp{Owen1977}), which is a generalization of the Shapley value suited to the case where features can be grouped into a priori coalitions. An example is given of machine learning-based decision models in automatic driving in which the nature of traffic naturally leads to the appearance of a partitioning of features in which those associated with the vehicle would form an a priori coalition.  
Table~\ref{tab:Literature_dep} provides a summary of references that propose and apply methods from cooperative game theory to enhance the interpretation of machine learning models in cases where feature dependence is present.

\section{Preliminaries}\label{sec:preli}
We start this section by introducing the fundamental concepts regarding datasets and influence measures necessary for the subsequent formal presentation of our model. Following this, we revisit key principles from cooperative game theory, given the close relationship between the influence measure proposed in this paper and a well-established solution concept within cooperative games.

\subsection{Datasets and influence measures}\label{subsec:prel_datasets}
Let $X = \{X_1, \dots, X_k\}$ be the set of features, with $K = \{1, \dots, k\}$ the set of indices of the features, and $Y$ a response variable. Also, let $\mathcal{A}_l$ denote the finite set of possible values or states that feature $X_l$, $l\in K$, can take, $\mathcal{A}=\prod_{l=1}^k \mathcal{A}_l$, and $\mathcal{B}$ the finite set of values that variable $Y$ can take. 
We will make use of datasets obtained from finite sets, denoted generally as $N$, consisting of $n$ individuals. That is, we have samples, in the form of $\mathcal{M} = \{(X^i, Y^i)\}_{i=1}^n$, where $X^i = (X^i_1, \dots, X_k^i)$ and $Y^i$, $i\in N$, are the observed values of the features and 
the prediction of the response variable, respectively, corresponding to individual $i$. These predictions are obtained by a classifier trained on the same set of individuals from which the true response variable values were observed. Note that $\cup_{i=1}^{n}{X^i}\subseteq \mathcal{A}$, so that the set of individuals can be identified with a set of different feature profiles. 

Formally, a dataset is a three-tuple $(X, Y, \mathcal{M})$ where $(X, Y)$ is a $k$-dimensional features vector and a response variable and $\mathcal{M}$ is a sample of size $n$. Given a dataset $(X, Y, \mathcal{M})$ with a binary response, i.e., $\abs{\mathcal{B}}=2$,  \cite{Datta2015} define a measure of influence of the features, for each $l\in K$, as follows:
\begin{equation*}
       \chi_l(X, Y, \mathcal{M}) = \sum_{(X^i, Y^i)\in \mathcal{M}} {\sum_{\substack{\bigl((X^i_{-l},\,a_{l}),\, b\bigl) \in \mathcal{M} \colon \\ a_{l}\in  \mathcal{A}_{l},\; b\in \mathcal{B} } }  }{|Y^i-b|}, 
\end{equation*}
where $(X^i_{-l},\,a_{l})=(X^i_1, \dots, X^i_{l-1}, a_l, X^i_{l+1}, \dots, X^i_k)$. Note that $\chi$ computes the number of times that a change in the state of feature $X_l$ causes a change in the response variable. 

%Next, we recall some related basic concepts of cooperative game theory, because the measure of influence introduced in this paper is related to a well-known solution concept of the so-called cooperative games.

\subsection{Cooperative games and values}\label{subsec:prel_games}
A cooperative game with transferable utility, often abbreviated as TU game, is a pair $(G, v)$ where $G$ denotes a finite set of players and $v$, the characteristic function, satisfies $v(\emptyset)=0$ and assigns a real number, $v(R)\in \mathbb{R}$, to each subset $R \subseteq G$. A TU game $(G, v)$ is called a simple game if: i) it is monotonous, meaning $v(R) \leq v(W)$ whenever $R\subseteq W\subseteq G$, ii) $v(R)\in \{0, 1\}$ for all $R\subseteq G$, and iii) $v(G)=1$. A simple game $(G, v)$ is a weighted majority game if there exists a vector of weights $w=(w_1, \dots, w_g)$ for the players, where $g=\abs{G}$, with $w_l \geq 0$, $1 \leq l\leq g$, and a positive real number $a \in \mathbb{R}^+$, referred to as the quota, such that $v(R) = 1$, $R\subseteq G$, if and only if $\sum_{l\in R}{w_l}\geq a$.

One of the focal aspects in cooperative game theory revolves around the definition of values. Some of these values serve as procedures for distributing the worth associated to the cooperation among players, as the Shapley value and the Owen value, while others can be used as ranking indices, as the Banzhaf value and the Banzhaf-Owen value. Given a TU game $(G, v)$ and $l\in G$, the Banzhaf value \citep{Banzhaf1965} is given by
 \begin{equation*}
     B_l(G, v)= \sum_{R\subseteq G\backslash \{l\}}\frac{1}{2^{g-1}}\cdot\Bigl(v(R\cup \{l\}) - v(R)\Bigl).
 \end{equation*}

\cite{Owen1981} extends the Banzhaf value to the class of TU games with a priori unions. A TU game with a priori unions \citep{Owen1977} is a three-tuple $(G, v, P)$, where $(G, v)$ is a TU game and $P= \{P_1, \dots, P_m\}$ is a partition of $G$ representing affinities among players, which might stem from familiar, political, or economic motives, among other factors. We denote by $M = \{1, \dots, m\}$ the set of indices of the coalitions, and we usually identify a coalition with its index.
Given a TU game with a priori unions $(G, v, P)$, $t\in M$, and $l\in P_t$, the Banzhaf-Owen value \citep{Owen1981} is given by
\begin{equation*}
BO_l(G, v, P)= \sum_{S\subseteq M\backslash \{t\}}
\sum_{R\subseteq P_t\backslash \{l\}}\frac{1}{2^{m-1}}\cdot\frac{1}{2^{|P_t|-1}}\cdot\Bigl(v(W\cup R\cup \{l\}) - v(W\cup R)\Bigl),
\end{equation*}
where $W = \cup_{u\in S}P_u$. Besides,  every coalition $W\cup R$, with $W = \cup_{u\in S}P_u$ and $R\subseteq P_t\setminus \{l\}$, is compatible with partition $P$ for $l$ in $P_t$.

\section{Main results}\label{sec:main}
In this section, we will consider $P = \{P_1, \dots, P_m\}$ a partition of $K$ that represents possible dependencies or interactions between the values of certain subsets of features. 
A trivial partition is $P^k=\{\{1\}, \dots,\{k\}\}$, where each subset of the partition is a singleton.

Now, a dataset with interactions is a four-tuple $(X, Y, P,\mathcal{M})$ where $(X, Y)$ is a $k$-dimen\-sional features vector and a response variable, $P$ is a partition of $K$, and $\mathcal{M}$ is a sample of size $n$; in other words, $(X, Y, \mathcal{M})$ is a dataset and $P$ is a partition of $K$, the set of indices of the features.
$D(X, Y)$ denotes the family of all datasets with interactions where $(X, Y)$ are the features vector and the response variable. Our goal is to make use of techniques used in classification problems to define a measure for studying the influence of features on the predicted value of the response variable under the assumption of possible dependencies between subsets of these features. First, let us state the formal definition of such an influence measure within this context. 

\begin{definition}
An influence measure for $D(X, Y)$ is a map, $I$, that assigns to every dataset with interactions, $(X, Y, P,\mathcal{M}) \in D(X, Y)$, a vector $I(X, Y, P,\mathcal{M})\in \mathbb{R}^k$. The measure $I_l(X, Y, P,\mathcal{M})\in \mathbb{R}$, $l\in K$, is a metric of the importance of $X_l$ in determining the predicted values of $Y$ over $\{(X^i)\}_{i=1}^n$.
\end{definition}

The main objective is to show that there is a unique measure of influence that satisfies a set of natural axioms, which we introduce and describe below. In what follows, we assume that $\abs{\mathcal{B}}=2$, i.e., the response variable $Y$ can only take two distinct values. 

A feature $X_l$, with $l\in K$, is said to be \textit{non-influential} in the dataset with interactions $(X, Y, P,\mathcal{M})$ if $Y^i=Y^j$ for all $i, j\in N$ such that $X^i_{-l}=X^j_{-l}$, where $X^i_{-l}$, with $i\in N$ and $l\in K$, denotes the vector $X^i$ after removing the $l$-th coordinate. First, we introduce the dummy property for influence measures.

\begin{itemize}
    \item[\textbf{(DP)}] {\bf Dummy property}. An influence measure $I$ satisfies the \textit{dummy property} if, for every $(X, Y, P,\mathcal{M})\in D(X, Y)$ and every non-influential feature in the dataset $(X, Y, P,\mathcal{M})$, $X_l$, with $l\in K$, it holds that 
$I_l(X, Y, P,\mathcal{M})=0$.
\end{itemize}

Given a dataset with interactions $(X, Y, P, \mathcal{M})$ and a bijective mapping $\sigma$ from $K$ to itself, %such that $\sigma(l)\in P_t$ for all $t\in M,\, l\in P_t$, 
we define $\sigma (X, Y, P, \mathcal{M})=(\sigma(X), Y, \sigma(P), \sigma (\mathcal{M}))$ in the natural way, consisting in relabelling the features according to $\sigma$, i.e., making the index of $l$, with $l\in K$, in the initial dataset become now $\sigma(l)$. We write $\sigma (\mathcal{M})=\{(\sigma(X^i), Y^i)\}_{i=1}^n$. Given a bijective mapping $\tau$ from $\mathcal{A}_l$, $l\in K$, to itself, we define
$\tau (X, Y, P, \mathcal{M})=(X, Y, P, \tau (\mathcal{M}))$ in a similar manner, consisting of relabelling the values of $\mathcal{A}_l$ according to $\tau$, i.e., making the value of $a_l$, with $a_l\in \mathcal{A}_l$, in the initial dataset becoming now $\tau(a_l)$. We write $\tau (\mathcal{M})=\{(\tau(X^i), Y^i)\}_{i=1}^n$.

These concepts enable us to introduce various notions of symmetry that measures of influence, like those examined here, should satisfy. 

\begin{itemize}
    \item[\textbf{(FSY)}] {\bf Feature symmetry}. An influence measure $I$ satisfies \textit{feature symmetry property} if, for every $(X, Y, P,\mathcal{M}) \in D(X, Y)$ such that $P=P^k$ and a bijective mapping $\sigma$ from $K$ to itself, 
%such that $\sigma(l)\in P_t$ for all $t\in M,\, l\in P_t$, 
it holds that $I_l(X, Y, P,\mathcal{M})=I_{\sigma(l)}(\sigma (X, Y, P,\mathcal{M}))$ for all $l\in K$.
\item[\textbf{(SSY)}] {\bf State symmetry}. An influence measure $I$ satisfies \textit{state symmetry property} if, for every $(X, Y, P,\mathcal{M})\in D(X, Y)$ such that $P=P^k$  and a bijective mapping $\tau$ from $A_l$, $l\in K$, to itself, 
%where $l\in P_t$, $t\in M$, 
it holds that $I_q(X, Y, P,\mathcal{M})  =I_q(\tau (X, Y, P,\mathcal{M}))$ for all $q\in K$.
\item[\textbf{(SY)}] {\bf Symmetry}. An influence measure $I$ satisfies \textit{symmetry property} if it satisfies both feature symmetry (FSY) and state symmetry (SSY).
\end{itemize}

Below we present some common properties that an influence measure should satisfy when considering a partition structure. Since $\abs{\mathcal{B}}=2$, we can assume that $\mathcal{B}=\{0,1\}$. Consequently, we define $W(\mathcal{M})=\{i\in N\, :\,  Y^i=1\}$ and $L(\mathcal{M})=\{i\in N\, :\, Y^i=0\}$ as the sets of sample profiles where the response variable takes the values 1 and 0, respectively. Thus, in general, given a set of individuals $N$, with $n\leq \abs{\mathcal{A}}$, $W, L\subseteq N$, and $W\cap L=\emptyset$, we can identify a sample with $(W,L)$ and a dataset with interactions with $(X,Y,P,(W,L))$.

\begin{itemize}
  \item[\textbf{(DU)}] {\bf Disjoint union}. An influence measure $I$ satisfies \textit{disjoint union property} if for every $(X, Y, P, (Q,R\cup R')), (X, Y, P, (R\cup R',Q))\in D(X,Y)$ where $X$ is a set of features and $Q, R$, and $R'$ are pairwise disjoint sets satisfying $\abs{Q\cup R\cup R'}\leq \abs{\mathcal{A}}$ for all $l\in K$, it holds that
\begin{equation*}
%\resizebox{0.82\textwidth}{!}{
    %$
    I_l(X, Y, P, (Q,R)) + I_l(X, Y, P, (Q,R')) = I_l(X, Y, P, (Q,R\cup R'))
    %$
%}
\end{equation*}
and
\begin{equation*}
%\resizebox{0.82\textwidth}{!}{
    %$
    I_l(X, Y, P, (R,Q)) + I_l(X, Y, P, (R',Q)) = I_l(X, Y, P, (R\cup R',Q)).
    %$
%}
\end{equation*}

\item[\textbf{(II)}] {\bf Indifference to interactions}. An influence measure $I$ satisfies \textit{indifference to interactions property} if for all $(X, Y, P,\mathcal{M})$ $\in D(X, Y)$, for all $t\in M$, with $l, q \in P_t$ and $l\neq q$, it holds that $I_l(X, Y, P,\mathcal{M}) = I_l(X, Y, P_{\_q},\mathcal{M})$, where $P_{\_q}$ is the partition that results from removing $q$ from its original subset to create a unitary subset, i.e., $P_{\_q}=\{P_1,\dots, P_{t-1}, P_t\backslash \{q\}, P_{t+1}, \dots, P_m, \{q\}\}.$ 
\item[\textbf{(RP)}] {\bf Relevance of dependent feature profiles}. An influence measure $I$ satisfies \textit{relevance of dependent feature profiles property} if for all 
$(X, Y, P, \mathcal{M})\in D(X, Y)$ and for all $t\in M$ such that $\abs{P_t}=1$, with $l\in P_t$, it holds that $I_l(X, Y, P, \mathcal{M}) = I_l(X, Y, P^k, \mathcal{M}^t)$, where 
\begin{equation}
\label{eq:emete}    
\mathcal{M}^t = \{(X^i, Y^i)\in \mathcal{M},\, i\in N \, :\, {\rm if}\, u\in M\backslash\{t\}\,, {\rm then}\, X^i_q \approx X^i_v\, {\rm for}\, {\rm all}\, q, v\in P_u\}
\end{equation}
is the subsample of $\mathcal{M}$ where features within the same union, except for $P_t$, precisely adhere to a pre-adjusted dependency model. In a union consisting of two binary features, the acceptable values for both features will be either identical or opposite, depending on whether the dependency is positive or negative.
\end{itemize}
 
Note that each of these properties extends specific axioms from the cooperative game theory literature to the context of influence measures. For instance, (DP) is an extension of the null player property utilized in the axiomatization of solutions for TU games, as evidenced in works such as \cite{feltkamp1995} for the Banzhaf value.  \cite{Datta2015} also employ this property from a binary classification perspective. Additionally, (FSY), (SSY), and (SY) share similarities with properties outlined in \cite{Datta2015}. These properties, in turn, expand upon those used in the axiomatization of the Banzhaf-Owen value for a TU game with a priori unions, among others, as seen in works like \cite{alonso2007b}. Furthermore, (DU) has the same essence as the union-intersection axiom employed in \cite{lehrer1988} to characterize the Banzhaf value for TU games. It also serves as a generalization of the axiom with the same name used in \cite{Datta2015}. Lastly, (II) and (RP) extend properties of indifference in unions and the quotient game for single-player unions introduced in \cite{alonso2007b} to axiomatize the Banzhaf-Owen value.

The following proposition extends the result of \cite{Datta2015} to the case of a set of categorical responses, when there is no partition structure on the affinities of the features.

\begin{proposition}\label{prop:datta}
An influence measure for $D(X, Y)$, $I$, satisfies (DP), (SY), and (DU) if and only if there exists a constant $C$ such that for every dataset with interactions $(X, Y, P^k, \mathcal{M})$ and every feature $l\in K$, 
\begin{equation}\label{eq:datta}
    I_l(X, Y, P^k, \mathcal{M})= C \cdot \sum_{(X^i, Y^i)\in \mathcal{M}}
    \, {\sum_{\substack{((X^i_{-l},\,a_{l}),\, b)\in \mathcal{M}\,:
    \\ a_{l}\in  \mathcal{A}_{l},\,\, b\in \mathcal{B}}}}|Y^i-b|.
\end{equation}
Moreover, it holds that $I_l(X, Y, P^k, \mathcal{M})=C \cdot\chi_l(X, Y, \mathcal{M})$.

\end{proposition}

Now, we present an axiomatic characterization of our influence measure in cases where a partition exists over the set of features, reflecting their dependencies and interactions.

\begin{remark}
It is important to note the analogy of this result with that obtained by \cite{alonso2007b} for the Banzhaf-Owen value in the context of TU games with a priori unions.
\end{remark}

\begin{theorem}[Existence and uniqueness]\label{th:1}
An influence measure for $D(X,Y)$, $I$, satisfies (DP), (SY), (DU), (II), and (RP) if and only if there exists a constant $C$ such that for every dataset with interactions $(X, Y, P, \mathcal{M})$, $t\in M$, $l\in P_t$, it holds that
\begin{equation}\label{eq:measure}
    I_l(X, Y, P, \mathcal{M})= C \cdot \Psi_l(X, Y, P, \mathcal{M}), 
\end{equation}
where
\begin{equation*}\label{eq:psi}
    \Psi_l(X, Y, P, \mathcal{M}) = \sum_{(X^i, Y^i)\in \mathcal{M}^t}
    \, {\sum_{\substack{
    ((X^i_{-l},\,a_{l}),\, b)\in \mathcal{M}^t\,:
    \\ a_{l}\in  \mathcal{A}_{l},\,\, b\in \mathcal{B}
    }
    }
    }|Y^i-b|.
\end{equation*}
%where $(X^i_{-l},\,a_{l})=(X^i_1, \cdots, X^i_{l-1}, a_l, X^i_{l+1}, \cdots, X^i_k)$.
\end{theorem}

\section{The influence measure and games with a priori unions}\label{sec:b-o}

In this section, we study the proposed measure of influence from a game-theoretical perspective. For this purpose, we mainly follow the ideas in \cite{Datta2015}.%, which starts by considering a TU game associated to any sample.

First, we will consider a TU game associated to any sample. %$\mathcal{M} = \{(X^i, Y^i)\}_{i=1}^n$. 
Let $\mathcal{M} = \{(X^i, Y^i)\}_{i=1}^n$ be a sample such that $\mathcal{A}=\cup_{i=1}^nX^i$ and $\abs{\mathcal{A}_l}=2$, for all $l\in K$. Thus, the sample $\mathcal{M}$ corresponds to the TU game $(K,v^{\mathcal{M}})$ defined, for all $R\subseteq K$, by
\begin{equation}\label{eq:game}
    v^{\mathcal{M}}(R)=Y^i \Longleftrightarrow \exists\, i\in N\; \textnormal{such that}\; 
    X^i_l= \left\{ \begin{array}{ll} 1 & {\rm if }\;  l\in R, \\ 0 &   {\rm otherwise}. \end{array} \right.
\end{equation}
If, moreover, we particularly consider the case where $\abs{\mathcal{B}}= 2$, we can assume that the characteristic function of the game fulfills $v^{\mathcal{M}}(R)\in \{0, 1\}$, for all $R \subseteq K$.

Using this proposal of TU game, \cite{Datta2015} innovatively relate their introduced influence measure to the Banzhaf value of the TU game $(K,v^{\mathcal{M}})$. More precisely, when the sample $\mathcal{M}$ corresponds to the TU game $(K,v^{\mathcal{M}})$ and $\abs{\mathcal{B}}= 2$, it follows that
\begin{equation}\label{eq:daba}
    B(K, v^{\mathcal{M}})=\dfrac{\chi(X, Y, \mathcal{M})}{\abs{\mathcal{A}}},
\end{equation}
indicating that the influence measure coincides with the (raw) Banzhaf value.

However, as previously justified, it is natural to assume the presence of a partition of the features describing potential dependencies among them. Consequently, a dataset with interactions $(X,Y,P,\mathcal{M})$, with $\mathcal{M}$ corresponding to a TU game $(K, v^{\mathcal{M}})$ and $P$ a partition of $K$, can be identified with a TU game with a priori unions $(K,v^{\mathcal{M}},P)$.

The following result extends Equation~\eqref{eq:daba} to the case of dependent features and shows that the  influence measure for potential dependent features introduced in this paper (see Equation~\eqref{eq:measure}) is a generalization of the Banzhaf-Owen value for simple games with a priori unions.

\begin{proposition}\label{prop:bo}
Let $(X, Y, P, \mathcal{M})$ be a dataset with interactions where the associations between the features within each union are positive. Suppose that $\mathcal{A}=\cup_{i=1}^nX^i$, $\abs{\mathcal{A}_l}=2$ for all $l\in K$, and $\abs{\mathcal{B}}=2$. If $(K, v^{\mathcal{M}})$ is the TU game corresponding to $\mathcal{M}$, then it holds that,
for every $t\in M$, $l\in P_t$, 
\begin{equation}\label{BO_psi}
    BO_l (K, v^{\mathcal{M}}, P)= \dfrac{\Psi_l(X, Y, P, \mathcal{M})}{\abs{\mathcal{A}^t}},
\end{equation}
where
$\mathcal{A}^t = \{a=(a_1,\dots,a_k)\in \mathcal{A}\, :\, {\rm if}\, u\in M\backslash\{t\}\; {\rm then}\; a_q = a_v\; {\rm for}\, {\rm all}\; q, v\in P_u\}$.
\end{proposition}

Below, we consider an example that illustrates the previous result.  

\begin{example}\label{ex:xogo}
Let us take the dataset with interactions $(X, Y, P, \mathcal{M})$ where $X=\{X_1, X_2, X_3, X_4\}$, that is, $k=4$, and such that
$\mathcal{A}_l = \mathcal{B} = \{0, 1\}$ for all $l\in K$. Suppose that $\mathcal{A}=\cup_{i=1}^nX^i$ and $n=16$. In addition, for $i\in N$, we have that
\begin{equation*}
    Y^i= \left\{ 
    \begin{array}{ll} 
        1 & \textrm{if }\; X^i\in \{0011, 1110, 1011, 0111, 1111\}, \\
        0 &  \textrm{otherwise.} 
    \end{array} 
    \right.
\end{equation*}
The dataset with interactions $(X, Y, P, \mathcal{M})$ can be identified with a weighted majority game with a priori unions $(K,v^{\mathcal{M}},P)$ where, for example, the vector of weights is $w$\,=\,(1, 1, 4, 3) and the quota is $a=6$. 
\begin{table}[h!]
\centering
\resizebox{12cm}{!}{
    \begin{tabular}{c|c|c} \hline 
    Scenario & Feature set partitions ($P$) & Influence measures of the features ($BO$)\\ \hline 
    1 & $\{\{1\},\{2\},\{3\},\{4\}\}$ & (0.125, 0.125, 0.625, 0.375) \\  
    2 &$\{\{1\},\{2\},\{3,4\}\}$ & (0.000, 0.000, 0,625, 0.375) \\
    3 &$\{\{1\},\{3\},\{2,4\}\}$ & (0.000, 0.125, 0.500, 0.375) \\
    4 &$\{\{1\},\{4\},\{2,3\}\}$ & (0.250, 0.125, 0.625, 0.250) \\
    5 &$\{\{2\},\{3\},\{1,4\}\}$ & (0.125, 0.000, 0.500, 0.375) \\
    6 &$\{\{2\},\{4\},\{1,3\}\}$ & (0.125, 0.250, 0.625, 0.250) \\
    7 &$\{\{3\},\{4\},\{1, 2\}\}$ & (0.125, 0.125, 0.750, 0.250) \\
    8 &$\{\{1, 2\}, \{3, 4\}\}$ & (0.000, 0.000, 0.750, 0.250)\\
    9 &$\{\{1, 3\}, \{2, 4\}\}$ & (0.000, 0.250, 0.500, 0.250) \\
    10 &$\{\{1, 4\}, \{2, 3\}\}$ & (0.250, 0.000, 0.500, 0.250) \\
    11 &$\{ \{1\},\{2, 3, 4\} \}$ & (0.000, 0.125, 0.625, 0.375) \\
    12 &$\{ \{2\},\{1, 3, 4\} \}$ & (0.125, 0.000, 0.625, 0.375) \\
    13 &$\{ \{3\},\{1, 2, 4\} \}$ & (0.125, 0.125, 0.500, 0.375) \\
    14 &$\{ \{4\},\{1, 2, 3\} \}$ & (0.125, 0.125, 0.625, 0.000) \\
    15 &$\{\{1, 2, 3, 4\}\}$ & (0.125, 0.125, 0.625, 0.375) \\ \hline
    \end{tabular}
}
\caption{Numerical results for our influence measure in the 15 possible scenarios of Example~\ref{ex:xogo}.}
\label{tab:xogo}
\end{table}
Table~\ref{tab:xogo} displays the influence measure introduced in this paper for the above-specified dataset, considering all possible partitions of the feature set. By  Proposition~\ref{prop:bo}, this influence measure coincides with the Banzhaf-Owen value. All results have been computed using the \texttt{R} library \texttt{powerindexR}\footnote{More information is available at \url{https://cran.r-project.org/web/packages/powerindexR/index.html}.}. 

From Table~\ref{tab:xogo}, we can clearly observe the behavior of various properties of our influence measure. Notably, in scenarios 1 and 2, features 3 and 4 exhibit identical influence, reflecting the property of indifference to interactions. This property is further confirmed by looking at the influence of features 3 and 4 in scenarios 2 and 12. Additionally, in scenario 1, where all feature unions have size 1, every feature is influential. However, when dependencies between groups of two or three features are considered, some features become non-influential, as demonstrated in scenarios 9 and 10.
\end{example}

In general, obtaining the exact Banzhaf values entails high computational complexity, specially for large sets of players, as the number of elements to be  evaluated increases exponentially (cf. \citealp{deng1994complexity}). While there are procedures to compute Banzhaf values for certain classes of games, which mitigate this issue, our approach does not benefit from such optimizations. Our influence measure requires the evaluation of all elements in the considered sample $\mathcal{M}$.

\section{Numerical results}\label{sec:results}
This section analyzes the performance of our influence measure on two different datasets that will be further described below. The proposed influence measure has been implemented in \texttt{R} {4.3.3} and a set of numerical experiments were run on a quad-core Intel i7-8665U CPU with 16 GB RAM.

We compute our influence measure for each dataset in several scenarios. On the one hand, we consider three different predictive models from the families of random forests (RFs), support vector machines (SVMs), and logistic regression (LR), chosen for their well-documented performance as classifier types \citep{fernandez2014}. We use the base implementations of these models provided by the \texttt{RWeka}\footnote{More information is available at \url{https://cran.r-project.org/web/packages/RWeka/RWeka.pdf}.} library in \texttt{R} software, specifically the \citeauthor{breiman2001}'s random forest classifier \citep{breiman2001}, the \citeauthor{Platt1998}'s sequential minimal optimization (SMO) algorithm for training a support vector machine classifier \citep{Platt1998}, and a modified \citeauthor{leCessie1992}'s multinomial logistic regression model with a ridge estimator \citep{leCessie1992}, respectively. 
On the other hand, we consider three distinct coalitional structures: 1) a partition formed by singletons, 2) a manually selected partition, and 3) a partition based on hierarchical clustering. For this latter case, the number of clusters is selected to match the number of coalitions in the manual partition, and we use the Jaccard distance measure, given the binary nature of our datasets.

The following subsections report the computational study. Subsection~\ref{subsec:cars} explores the influence of certain factors on the severity of a car crash, while Subsection~\ref{subsec:spotify} examines the influence of listening to various music groups or singers on the likelihood of listening to other artists. We will also discuss how the chosen classifiers and partitions impact feature rankings.
All the datasets used are available at \url{https://github.com/LauraDavilaPena/GT-based_IM}.

\subsection{On the analysis of car crash fatalities}\label{subsec:cars}

In this section, we apply our proposed influence measure to the relevant context of occupant safety in car crashes. The severity of an accident can vary significantly based on the type of collision. For instance, a head-on collision differs substantially from a side-impact collision or a vehicle rollover. Furthermore, the risk of fatality in a car older than ten years old is twice as high as in a newer vehicle. We will specifically look at the variables that influence the likelihood of fatalities in vehicle accidents. 

\begin{table}[ht!]
\centering \resizebox{0.88\textwidth}{!}{
    \begin{tabular}{p{0.15cm}p{1.5cm}p{12.5cm}}
        \noalign{\smallskip}\hline\noalign{\smallskip}
      \multicolumn{2}{l}{Response} & Description \\
      \noalign{\smallskip}\hline\noalign{\smallskip}
     & \texttt{deceased} &  Binary variable indicating whether the person involved in the car accident is deceased (1) or not (0). \\
     \noalign{\smallskip}\hline\noalign{\smallskip}
     \multicolumn{2}{l}{Feature} & Description \\
     \noalign{\smallskip}\hline\noalign{\smallskip}
    1 & \texttt{dvcat} & Binary variable indicating whether the vehicle, at the moment of the accident, was traveling at a speed higher than 55 km/h (1) or not (0).\\[0.2em]
    2 & \texttt{airbag} & 
    Binary variable indicating whether the vehicle had an airbag system (1) or not (0). \\[0.2em]
    3 & \texttt{seatbelt} & 
    Binary variable indicating whether the person involved was wearing a seat belt (1) or not (0). \\[0.2em]
    4 & \texttt{frontal} & 
    Binary variable indicating whether the vehicle crash was frontal (1) or non-frontal (0). \\[0.2em]
    5 & \texttt{sex} & Binary variable indicating the sex of the person involved: 1 for male and 0 for female. \\[0.2em]
    6 & \texttt{abcat} & Binary variable indicating airbag activation: 1 if one or more airbags in the vehicle were activated (even if not deployed) and 0 if none were deployed (either due to malfunction or being disabled).\\[0.2em]
    7 & \texttt{occRole} & Binary variable indicating whether the person involved was the driver (1) or a passenger (0) of the vehicle.\\[0.2em]
    8 & \texttt{deploy} & Binary variable indicating whether the airbag functioned correctly (1) or was unavailable or not functioning (0). \\[0.2em]
    9 & \texttt{ageOFocc} &  Binary variable indicating whether the person was 30 years old or less (1) or over 30 years old (0).\\[0.2em]
    10 & \texttt{age} & Binary variable indicating whether the vehicle was 10 years old or more (1) or less than 10 years old (0).\\
    \noalign{\smallskip}\hline\noalign{\smallskip}
    \end{tabular}
    }
    \caption{Summary of the considered features in the analysis of car crash fatalities.}
    \label{tab:variablecars}
\end{table}

To achieve this, we examine the influence of various factors that help us describe the nature of a road accident. Specifically, we consider the variables included in the \texttt{nassCDS} dataset from \texttt{R} software. This dataset contains information on car crashes in the US in the period 1997-2002, reported by the police, in which there is an injury (to person or property) and at least one vehicle is towed. The data is limited to front-seat occupants and includes only a subset of the recorded variables, with additional restrictions. For our application, we have conveniently adapted this database to focus on binary variables. 
The transformed database consists of a sample of 17,565 observations for the 10 characteristics listed in Table~\ref{tab:variablecars}.

From the sample, we obtain the influence measure both with and without a coalitional structure of the features across various scenarios. Initially, we consider the case of partitioning with unitary elements, meaning we do not account for potential affinity relations between the features; this corresponds to the original case presented by \cite{Datta2015}. Next, we examine a constructed-by-design partition (CDP) where features 1, 4, and 9 are grouped in one block, 2, 6, and 8 in another block, and features 3, 5, 7, and 10 act individually. Finally, we adopt a coalitional structure specified by hierarchical clustering, fixing the number of clusters at six. This hierarchical clustering partition (HCP) joins features 2, 7, and 9 in the same block, features 4, 5, and 8 in another block, and leaves 1, 3, 6, and 10 alone. 
For each scenario, we compute the resulting influence measure (IM) using the above-specified random forest (RF), support vector machine (SVM), and logistic regression (LR) classifiers as predictive models. Table~\ref{tab:cars_measures} displays the numerical results.

\begin{table}[ht!]
    \centering\resizebox{0.88\textwidth}{!}{
    \begin{tabular}{p{1.3 cm}|p{1.3 cm}p{1.3 cm}p{1.3 cm}|p{1.3 cm}p{1.3 cm}p{1.3 cm}|p{1.3 cm}p{1.3 cm}p{1.3 cm}}\hline
    & \multicolumn{3}{c|}{\citeauthor{Datta2015}'s IM, \eqref{eq:datta}} &\multicolumn{3}{c|}{IM with CDP, \eqref{eq:measure}} & \multicolumn{3}{c}{IM with HCP, \eqref{eq:measure}}\\
    \hline
   Feature    	& RF	& SVM &	LR			& RF & SVM	&	LR	& RF & SVM	&	LR \\\hline
1	& 0.42801	& 0 &	0.19516	 &	0.63830	& 0 &	0.53216	&	0.03189	& 0 &	0.03826\\
2	& 0.00865	& 0 &	0.00159	&		0	&0& 	0		&	0.00798	& 0 &	0.00059\\
3	& 0.02847	& 0&         0.08289	&		0	& 0 &	0		&	0.01934 & 0	&	0.03491\\
4	& 0.11967	& 0 &	0.34728	&	0.15016	& 0 &	0.27097	&	0.13540 & 0	&	0.13117\\
5	& 0.03450	& 0 &	0.01697	&		0 & 0	&	0		&	0.00106 & 0 &	0\\
6	& 0.03507 & 0	&	0.00649	&	0.01045	& 0 &	0		&	0.00319 & 0	&	0\\
7	& 0		& 0 &      0		&  0  & 0 & 0   	&	0		&	0	&	0\\
8	& 0.03632 & 0	&	0.00330	&	0.03947 & 0	&	0	&	0.02977	& 0 &	0\\
9	& 0.02812 & 0	&	0.004787	&	0.05136 & 0	&	0.01663	&	0.02157 & 0	&	0.00177\\
10	& 0.01754	& 0 &	0.00137	&	0	& 0 & 0	&	0		&	0	&	0\\\hline
    \end{tabular}}
    \caption{Numerical results in the analysis of car crash fatalities.}
    \label{tab:cars_measures}
\end{table}

Based on these findings, several conclusions can be drawn. Features 1 and 4 consistently hold the top two positions in the ranking of the most influential features. Specifically, they refer to the vehicle's speed and whether the collision is frontal. Factors such as airbag functioning (feature 8), seat belt usage (feature 3), or the person's age (feature 9) rank third in the different rankings obtained, as can be seen in Table~\ref{tab:cars_measuresr}. Notably, regardless of the scenario examined, the individual's role (passenger or driver), represented by feature 7, appears to have no influence.

\begin{table}[ht!]
    \centering\resizebox{0.9\textwidth}{!}{
    \begin{tabular}{p{1.3 cm}|p{1.3 cm}p{1.3 cm}p{1.3 cm}|p{1.3 cm}p{1.3 cm}p{1.3 cm}|p{1.3 cm}p{1.3 cm}p{1.3 cm}}\hline
    & \multicolumn{3}{c|}{\citeauthor{Datta2015}'s IM, \eqref{eq:datta}} &\multicolumn{3}{c|}{IM with CDP, \eqref{eq:measure}} & \multicolumn{3}{c}{IM with HCP, \eqref{eq:measure}}\\
    \hline
   Position    	& RF	& SVM &	LR			& RF & SVM	&	LR	& RF & SVM	&	(LR)\\\hline
1	&1 & -& 4 & 1 &-& 1& 4&-&4\\
2	&4 & -& 1 & 4 &-& 4& 1&-&1\\
3	&8 & -& 3 & 9 &-& 9& 8&-&3\\\hline
    \end{tabular}}
    \caption{Top three features in the analysis of car crash fatalities.}
    \label{tab:cars_measuresr}
\end{table}

It is also worth noting that the SVM classifier does not yield conclusive results when using measures based on \citeauthor{Datta2015}'s methodology. In such cases, influence measures always return a value of 0. As discussed in the literature, the imbalanced distribution between 0s and 1s in the response variable can lead to influence measure outcomes that may not accurately reflect reality. As shown in Table~\ref{tab:cars_count}, SVM does not predict any occurrences of the value 1 for the response variable in this particular case study. 

\begin{table}[ht!]
    \centering\resizebox{0.5\textwidth}{!}{
    \begin{tabular}{p{1.3 cm}|p{1.3 cm}p{1.3 cm}p{1.3 cm}p{1.3 cm}}\hline
   Value    & \multicolumn{1}{l}{Original response}	& \multicolumn{1}{l}{RF}	& \multicolumn{1}{l}{SVM} &	\multicolumn{1}{l}{LR}	\\\hline
0 &  \multicolumn{1}{r}{16,788}	&  \multicolumn{1}{r}{17,496}	&  \multicolumn{1}{r}{17,565} &	 \multicolumn{1}{r}{17,464}	 \\
1 & \multicolumn{1}{r}{777}	& \multicolumn{1}{r}{69}	& \multicolumn{1}{r}{0} &	\multicolumn{1}{r}{101}	\\\hline
    \end{tabular}}
    \caption{Summary of the responses from both the original database and those predicted by the classifiers.}
    \label{tab:cars_count}
\end{table}

To check this conjecture, we draw a subsample from the original database consisting of 777 instances with a response value of 1 and 777 instances with a response value of 0 (balanced in response). Table~\ref{tab:cars_subsample} presents the numerical results for this ad hoc subsample. Interestingly, in this case, we observe that the influence measures based on \citeauthor{Datta2015}'s methodology are not zero for SVM's scenarios.

\begin{table}[ht!]
\centering\resizebox{0.9\textwidth}{!}{
    \begin{tabular}{p{1.3 cm}|p{1.3 cm}p{1.3 cm}p{1.3 cm}|p{1.3 cm}p{1.3 cm}p{1.3 cm}|p{1.3 cm}p{1.3 cm}p{1.3 cm}}\hline
    & \multicolumn{3}{c|}{\citeauthor{Datta2015}'s IM, \eqref{eq:datta}} &\multicolumn{3}{c|}{IM with CDP, \eqref{eq:measure}} & \multicolumn{3}{c}{IM with HCP, \eqref{eq:measure}}\\
    \hline
   Feature    	& RF	& SVM &	LR			& RF & SVM	&	LR	& RF & SVM	&	LR\\\hline
1 & 1.94852 & 2.01931 & 2.23681 & 2.70455 & 2.70455 & 4.04546 & 4.09669 & 4.09669 & 4.95165\\
2 & 0.22780 & 0.17503 & 1.24324 & 0.39351 & 0.49087 & 3.61460 & 0.34128 & 0.02202 & 0.45505\\
3 & 8.15187 & 3.60489 & 7.11454 & 6.87500 & 1.50000 & 0 & 15.96947 & 6.52417 & 14.97201\\
4 & 3.33333 & 6.39768 & 5.10296 
& 2.23864 & 1.97159 & 1.86932 & 4.59246 & 8.93836 & 7.13340\\
5 & 0.79022 & 1.10039 & 0 & 1.64286 & 0 & 0 & 0.48574 & 1.57314 & 0\\
6 & 3.78765 & 2.50708 & 4.77992 
& 1.08316 & 3.22110 & 3.78093 
& 9.48092 & 4.04580 & 10.22392\\
7 & 0 & 0 & 0 & 0 & 0 & 0 & 0 & 0 & 0\\
8 & 1.57143 & 2.04376 & 0 & 3.33874 & 3.39148 & 0 & 1.19411 & 1.67065 & 0\\
9 & 0.00257 & 0 & 0.01030 & 0.01136 & 0 & 0.04546 & 0 & 0 & 0.02202\\
10 & 0.71429 & 0.24710 & 0 & 0.19643 & 0 & 0 & 1.52163 & 0 & 0 \\\hline
    \end{tabular}}
    \caption{Numerical results in the analysis of car crash fatalities over the considered subsample.}
    \label{tab:cars_subsample}
\end{table}

A consistent finding shared with the original case is the lack of influence of feature 7  (the passenger's role) across all studied scenarios. Furthermore, Table~\ref{tab:cars_measuresr2} presents additional insights, showcasing the top three features.

\begin{table}[ht!]
    \centering\resizebox{0.9\textwidth}{!}{
    \begin{tabular}{p{1.3 cm}|p{1.3 cm}p{1.3 cm}p{1.3 cm}|p{1.3 cm}p{1.3 cm}p{1.3 cm}|p{1.3 cm}p{1.3 cm}p{1.3 cm}}\hline
    & \multicolumn{3}{c|}{\citeauthor{Datta2015}'s IM, \eqref{eq:datta}} &\multicolumn{3}{c|}{IM with CDP, \eqref{eq:measure}} & \multicolumn{3}{c}{IM with HCP, \eqref{eq:measure}}\\
    \hline
   Position    	& RF	& SVM &	LR			& RF & SVM	&	LR	& RF & SVM	&	LR\\\hline
1	& 3 & 4 & 3 & 3 & 8 & 1 & 3 & 4 & 3\\
2	& 6 & 3 & 4 & 8 & 1 & 6 & 6 & 3 & 6\\
3	& 4 & 6 & 6 & 1 & 6 & 2 & 4 & 1 & 4\\\hline
    \end{tabular}}
    \caption{Top three features in the analysis of car crash fatalities over the considered subsample.}
    \label{tab:cars_measuresr2}
\end{table}

Features 3, 4, and 6, which relate to seat belt usage, crash type, and airbag activation, respectively, are the most influential factors (with slight variations in order) for both the influence measure without coalitional structure or with coalitional structure provided by hierarchical clustering. It is reasonable to assume that safety features such as seat belts and airbags play a pivotal role in ensuring passenger safety. In the hierarchical clustering scenario with SVM, feature 6 is replaced by feature 1. More disparities are found when considering the partition constructed by mere feature observation. Only feature 1 ranks in the top three across all scenarios. Feature 8 emerges when RF and SVM are selected, while feature 6 appears alongside SVM and LR. Feature 3 completes the ranking given by RF, and feature 2 does so for LR.

As a final thought, the numerical results presented here seem to be strongly dependent on the selected sample and the methodology employed in each case. Nevertheless, it is notable that the ultimate discrepancies are minimal in terms of rankings, with only a few variations in the positions. The incorporation of information regarding feature affinities through partitions may justify this fact.

\subsection{On the analysis of musical taste in Spotify}\label{subsec:spotify}

Next, we use our proposed influence measure to analyze musical tastes from a Spotify database. Our goal is to identify the most influential singers or music groups. The original dataset contains a total of 285 artists and 1226 users and includes which users have listened to which artists. For the purposes of our study, we have decided to consider only the 75 most listened-to artists, all of them with an audience share of no less than 5\%, to ensure results in a reasonable computational time and to avoid examining those artists who are already of little interest.

\begin{table}[h!]
    \centering\centering\resizebox{0.9\textwidth}{!}{\begin{tabular}{l|cc||l|cc}\hline
Artist	&	Decade	& \% of listeners	&		Artist	&	Decade	&		\% of listeners\\\hline
\texttt{linkin.park}	&	1990	&	16.2	&		\texttt{johnny.cash}	&	1950	&	7.3	\\
\texttt{coldplay}	&	1990	&	16.1	&		\texttt{kings.of.leon}	&	2000	&	7.3	\\
\texttt{red.hot.chili.peppers}	&	1980	&	15.1	&		\texttt{amy.winehouse}	&	2000	&	7.2	\\
\texttt{rammstein}	&	1990	&	15\phantom{.0}	&		\texttt{depeche.mode}	&	1980	&	7.2	\\
\texttt{system.of.a.down}	&	1990	&	13.2	&		\texttt{bullet.for.my.valentine}	&	1990	&	7\phantom{.0}	\\
\texttt{metallica}	&	1980	&	12.2	&		\texttt{in.extremo}	&	1990	&	7\phantom{.0}	\\
\texttt{die.toten.hosen}	&	1980	&	11.8	&		\texttt{blink.182}	&	1990	&	6.9	\\
\texttt{billy.talent}	&	2000	&	10.8	&		\texttt{slipknot}	&	1990	&	6.7	\\
\texttt{the.killers}	&	2000	&	10.8	&		\texttt{death.cab.for.cutie}	&	1990	&	6.5	\\
\texttt{the.beatles}	&	1960	&	10.6	&		\texttt{daft.punk}	&	1990	&	6.4	\\
\texttt{jack.johnson}	&	2000	&	10\phantom{.0}	&		\texttt{limp.bizkit}	&	1990	&	6.3	\\
\texttt{muse}	&	1990	&	10\phantom{.0}	&		\texttt{sum.41}	&	2000	&	6.2	\\
\texttt{beatsteaks}	&	1990	&	9.5	&		\texttt{fall.out.boy}	&	2000	&	6.1	\\
\texttt{foo.fighters}	&	1990	&	9.5	&		\texttt{schandmaul}	&	1990	&	6.1	\\
\texttt{nirvana}	&	1980	&	9.5	&		\texttt{kanye.west}	&	2000	&	6\phantom{.0}	\\
\texttt{radiohead}	&	1990	&	9.4	&		\texttt{seeed}	&	1990	&	6\phantom{.0}	\\
\texttt{arctic.monkeys}	&	2000	&	9.1	&		\texttt{tenacious.d}	&	1990	&	6\phantom{.0}	\\
\texttt{placebo}	&	1990	&	9.1	&		\texttt{rihanna}	&	2000	&	6\phantom{.0}	\\
\texttt{bloc.party}	&	2000	&	8.8	&		\texttt{papa.roach}	&	1990	&	5.8	\\
\texttt{evanescence}	&	1990	&	8.8	&		\texttt{portishead}	&	1990	&	5.8	\\
\texttt{rise.against}	&	2000	&	8.7	&		\texttt{rage}	&	1990	&	5.8	\\
\texttt{the.kooks}	&	2000	&	8.7	&		\texttt{franz.ferdinand}	&	2000	&	5.7	\\
\texttt{mando.diao}	&	2000	&	8.6	&		\texttt{marilyn.manson}	&	1980	&	5.6	\\
\texttt{the.white.stripes}	&	1990	&	8.4	&		\texttt{nelly.furtado}	&	2000	&	5.6	\\
\texttt{deichkind}	&	1990	&	8.3	&		\texttt{queen}	&	1970	&	5.5	\\
\texttt{incubus}	&	1990	&	8.3	&		\texttt{bob.marley}	&	1960	&	5.5	\\
\texttt{farin.urlaub}	&	1990	&	8.2	&		\texttt{feist}	&	1990	&	5.5	\\
\texttt{in.flames}	&	1990	&	8.2	&		\texttt{massive.attack}	&	1980	&	5.4	\\
\texttt{clueso}	&	2000	&	8.1	&		\texttt{queens.of.the.stone.age}	&	1990	&	5.4	\\
\texttt{peter.fox}	&	2000	&	8\phantom{.0}	&		\texttt{iron.maiden}	&	1970	&	5.2	\\
\texttt{the.offspring}	&	1980	&	8\phantom{.0}	&		\texttt{avril.lavigne}	&	2000	&	5.1	\\
\texttt{air}	&	1990	&	7.8	&		\texttt{amon.amarth}	&	1990	&	5.1	\\
\texttt{subway.to.sally}	&	1990	&	7.7	&		\texttt{apocalyptica}	&	1990	&	5.1	\\
\texttt{nightwish}	&	1990	&	7.7	&		\texttt{gorillaz}	&	1990	&	5.1	\\
\texttt{the.prodigy}	&	1990	&	7.4	&		\texttt{nine.inch.nails}	&	1980	&	5.1	\\
\texttt{disturbed}	&	1990	&	7.3	&		\texttt{oasis}	&	1990	&	5.1	\\
\texttt{ac.dc}	&	1970	&	7.3	&		\texttt{children.of.bodom}	&	1990	&	5\phantom{.0}	\\
\texttt{green.day}	&	1980	&	7.3	&			&		&		\\\hline
    \end{tabular}}
    \caption{List of the 75 most listened-to music bands in the database, with the prescribed partition by decade of the band's creation and percentage of users who listen to it.}
    \label{tab:list_music_bands}
\end{table}

Table~\ref{tab:list_music_bands} presents the selected artists, showcasing the percentage of listeners and the decade in which their band was created. This will constitute the basis for our CDP, which groups artists based on the decade of their foundation. %We then have $P^{d}=\cup_{k=1}^{6}{P^{d}_{k}}$, where $P^{d}_{k},\: k=1,\dots,6$, consists of all artists whose band was created in 1950, 1960, 1970, 1980, 1990, and 2000, respectively.
Table~\ref{tab:music_bands_decade} summarizes the number of artists in each group. 
%We do not care when the user listened to the song, as we assume that music tastes are fairly static over time.

\begin{table}[h!]
    \centering
    \begin{tabular}{c|cccccc}\hline
Decade & 1950 & 1960 & 1970 & 1980 & 1990 & 2000 \\ \hline
\# of artists & 1 & 2 & 3 &10& 40& 19\\ \hline
    \end{tabular}
    \caption{Distribution of artists per decade.}
    \label{tab:music_bands_decade}
\end{table}

To apply our methodology, we first need to define our classification problems. We take  10 of these artists as the response variable for each of our problems, and our objective will be to study whether or not listening to the other 74 artists affects this target.  In this example, we compare the results obtained in two different scenarios. First, we consider the case where the partition is determined by the decade of the band's foundation, the CDP. Then, as in Subsection~\ref{subsec:cars}, we consider hierarchical clustering techniques to obtain a partition of the bands for each of the considered response variables, that is, the HCP. Note that the large number of bands involved makes it computationally infeasible to calculate exactly the influence measure when the partition is composed of singletons, that is, \citeauthor{Datta2015}'s influence measure \eqref{eq:datta}.

\begin{table}[h!]
    \centering\resizebox{0.88\textwidth}{!}{
    \begin{tabular}{p{2.55cm}|c|p{2.8 cm}p{3 cm}p{3 cm}|p{3 cm}p{3 cm}p{2.75 cm}}\hline
    \multicolumn{2}{c}{ } & \multicolumn{3}{|c|}{IM with CDP, \eqref{eq:measure}}& \multicolumn{3}{c}{IM with HCP, \eqref{eq:measure}}\\\hline
\multicolumn{2}{c|}{ }  & Position 1 & Position 2 & Position 3  & Position 1 & Position 2 & Position 3 \\\hline
\texttt{coldplay} & RF & \texttt{avril.lavigne} & \texttt{muse} & \texttt{air} & \texttt{avril.lavigne} & \texttt{muse} & \texttt{air}\\
& SVM & - & - & - & - & - & - \\
& LR & - & - & - & - & - & - \\\hline
\texttt{metallica} & RF & \texttt{nightwish} & \texttt{rage} & \texttt{evanescence} & \texttt{nightwish} & \texttt{iron.maiden} & \texttt{amon.amarth}\\
& SVM & - & - & - & - & - & - \\
& LR & - & - & - & \texttt{nightwish} & \texttt{iron.maiden} & - \\\hline
\texttt{the.killers} & RF & \texttt{green.day} & \texttt{the.white.stripes} & \texttt{kings.of.leon} & \texttt{green.day} & \texttt{franz.ferdinand} & \texttt{evanescence}\\
& SVM & - & - & - & - & - & - \\
& LR & - & - & - & - & - & -  \\\hline
\texttt{the.beatles} & RF & \texttt{kings.of.leon} & \texttt{air} & \texttt{franz.ferdinand} & \texttt{kings.of.leon} & \texttt{muse} & \texttt{apocalyptica}\\
& SVM & - & - & - & - & - & -  \\
& LR & - & - & - & - & - & -  \\\hline
\texttt{muse} & RF & \texttt{avril.lavigne} & \texttt{evanescence} & \texttt{placebo} & \texttt{evanescence} & \texttt{avril.lavigne} & - \\
& SVM & - & - & - & - & - & -  \\
& LR & - & - & - & - & - & -  \\\hline
\texttt{ac.dc} & RF & \texttt{rage} & \texttt{air} & \texttt{the.prodigy} & \texttt{rage} & \texttt{evanescence} & \texttt{air}\\
& SVM & - & - & - & - & - & -  \\
& LR & - & - & - & - & - & - \\\hline
\texttt{amy.winehouse} & RF & \texttt{kings.of.leon} & \texttt{nelly.furtado} & \texttt{peter.fox} & \texttt{kings.of.leon} & \texttt{bob.marley} & \texttt{the.beatles} \\
& SVM & - & - & - & - & - & -  \\
& LR & - & - & - & - & - & -  \\\hline
\texttt{rihanna} & RF & \texttt{limp.bizkit} & \texttt{the.prodigy} & \texttt{amy.winehouse} & \texttt{evanescence} & \texttt{seeed} & \texttt{nelly.furtado} \\
& SVM & \texttt{nelly.furtado} & \texttt{kanye.west} & - & \texttt{seeed} & \texttt{nelly.furtado} & \texttt{kanye.west}\\
& LR & \texttt{nelly.furtado} & \texttt{kanye.west} & - & \texttt{seeed} & \texttt{nelly.furtado} & \texttt{kanye.west} \\\hline
\texttt{queen} & RF & \texttt{amon.amarth} & \texttt{the.prodigy} & \texttt{muse} & \texttt{feist} & \texttt{muse} & \texttt{the.beatles}\\
& SVM & - & - & - & - & - & - \\
& LR & - & - & - & - & - & - \\\hline
\texttt{bob.marley} & RF & \texttt{seeed} & \texttt{jack.johnson} & \texttt{the.killers} & \texttt{seeed} & \texttt{the.killers} & \texttt{peter.fox} \\
& SVM & - & - & - & - & - & - \\
& LR & - & - & - & - & - & -  \\\hline
    \end{tabular}}
    \caption{Top three bands in the analysis of musical taste in Spotify.}
    \label{tab:ranking_decades}
\end{table} 

Table~\ref{tab:ranking_decades} summarizes the three most influential bands for each scenario. The numerical results show some similarities, which are detailed below. For instance, according to RF, Avril Lavigne, Muse, and Air are the three most influential bands for listening to Coldplay. In the case of The Killers, Green Day is the most influential band; for The Beatles, Kings of Leon; for Muse, Avril Lavigne and Evanescence swap positions when using the CDP and the HCP; for AC/DC, Rage and Air are both in the top three; for Amy Winehouse, the most influential band is Kings of Leon; and for Queen, the only coincidence is that Muse is in position 3 when using the CDP and in position 2 when using the HCP. Finally, Seeed is the most influential band for listening to Bob Marley, although The Killers are in positions 3 and 2 for the CDP and the HCP, respectively. In all these cases, no other bands are influential for SVM and LR. However, for Metallica, Nightwish is the most influential band according to RF and both partition structures, as well as LR with the HCP. Also in this latter case, Iron Maiden is in the second position under RF and LR. When analyzing the case of Rihanna, Nelly Furtado ranks first for the CDP with SVM and LR, second for the HCP with SVM or LR, and third for the HCP with RF. Seeed is the most influential when using the HCP with SVM and LR, and ranks second with RF. We also highlight the case of Kanye West, in positions 2 and 3 for both SVM and LR with CDP and HCP, respectively.

\subsubsection{The case of the 15 most listened-to music bands}\label{subs:15artists}

As the reader can check in the previous section, the SVM and LR classifiers also did not yield conclusive results in the analysis of musical tastes. Therefore, we have narrowed our study to focus only on the 15 most listened-to bands in our original dataset. According to Table~\ref{tab:list_music_bands}, these bands have been listened to by more than 9.5\% of the users included in the database. 

The outline of this section follows our previous approach. We once again apply  RF, SVM, and LR as classifiers on the dataset under consideration, and we obtain the influence measure using the two aforementioned partitions: CDP and HCP. Despite focusing on the 15 most listened-to bands, we encountered inconclusive results with SVM when assessing influences.

First, we analyze the case of using the CDP. The corresponding numerical results are shown in Table~\ref{tab:musicdecades}. For each classification problem (by columns), we highlight the most influential band in bold. Using RF, the following conclusions can be drawn. Linkin Park is the most influential band for listening to System of a Down, Metallica and Die Toten Hosen; Rammstein, to Billy Talent; System of a Down, to Rammstein; Die Toten Hosen, to Nirvana; Billy Talent, to Linkin Park and Muse; Muse, to Beatsteaks; Beatsteaks, to Red Hot Chilli Peppers, The Killers, and Jack Johnson; Foo Fighters, to Coldplay; and Nirvana, to The Beatles. When considering LR, Coldplay and Foo Fighters are both the most influential bands on Beatsteaks; Rammstein, on Metallica and Die Toten Hosen; Metallica, on The Beatles; Muse, on Rammstein; Beatsteaks, on The Killers and Foo Fighters; and Foo Fighters, on Coldplay, Billy Talent, Muse, and Beatsteaks. Only Foo Fighters and Beatsteaks share a bilateral relation where each band is mutually the most influential on the other.

\begin{table}[h!]
    \centering\resizebox{0.9\textwidth}{!}{
    \begin{tabular}{p{4.2 cm}|p{1.1 cm}p{1.1 cm}p{1.1 cm}p{1.1 cm}p{1.1 cm}p{1.1 cm}p{1.1 cm}p{1.1 cm}p{1.1 cm}p{1.1 cm}p{1.1 cm}p{1.1 cm}p{1.1 cm}p{1.1 cm}p{1.1 cm}}\hline
   RF & \rotatebox[origin=c]{90}{\texttt{linkin.park}} & \rotatebox[origin=c]{90}{\texttt{coldplay}} & \rotatebox[origin=c]{90}{\mbox{ }\texttt{red.hot.chili.peppers}\mbox{ }} & \rotatebox[origin=c]{90}{\texttt{rammstein}} & \rotatebox[origin=c]{90}{\texttt{system.of.a.down}} & \rotatebox[origin=c]{90}{\texttt{metallica}} & \rotatebox[origin=c]{90}{\texttt{die.toten.hosen}} & \rotatebox[origin=c]{90}{\texttt{billy.talent}} & \rotatebox[origin=c]{90}{\texttt{the.killers}} & \rotatebox[origin=c]{90}{\texttt{the.beatles}} & \rotatebox[origin=c]{90}{\texttt{jack.johnson}} & \rotatebox[origin=c]{90}{\texttt{muse}} & \rotatebox[origin=c]{90}{\texttt{beatsteaks}} & \rotatebox[origin=c]{90}{\texttt{foo.fighters}} & \rotatebox[origin=c]{90}{\texttt{nirvana}}\\\hline
\texttt{linkin.park} & - & 0.0000 & 0.0028 & 0.5128 &\textbf{1.0256} &\textbf{0.1114} &\textbf{0.8350} & 0.0282 & 0.0000 & 0.0000 & 0.0028 & 0.0090 & 0.0271 & 0.0000 & 0.0000\\
\texttt{coldplay} & 0.0000 & - & 0.0693 & 0.0241 & 0.0271 & 0.0000 & 0.0226 & 0.0367 & 0.0905 & 0.0000 & 0.0056 & 0.0090 & 0.0211 & 0.0754 & 0.0000\\
\texttt{red.hot.chili.peppers} & 0.0000 & 0.0000 & - & 0.3110 & 0.0000 & 0.0000 & 0.0078 & 0.0038 & 0.0000 & 0.0117 & 0.0109 & 0.0113 & 0.0000 & 0.0000 & 0.0313\\
\texttt{rammstein} & 0.0000 & 0.0030 & 0.0166 & - & 0.7903 & 0.0167 & 0.6996 & \textbf{0.0423} & 0.0082 & 0.0000 & 0.0056 & 0.0000 & 0.0030 & 0.0030 & 0.0057\\
\texttt{system.of.a.down} & 0.0000 & 0.0121 & 0.0111 & \textbf{0.7149} & - & 0.0251 & 0.1128 & 0.0000 & 0.0137 & 0.0000 & 0.0084 & 0.0000 & 0.0121 & 0.0060 & 0.0115\\
\texttt{metallica} & 0.0000 & 0.0000 & 0.0156 & 0.4205 & 0.0000 & - & 0.0078 & 0.0000 & 0.0000 & 0.0117 & 0.0219 & 0.0226 & 0.0000 & 0.0000 & 0.0234\\
\texttt{die.toten.hosen} & 0.0000 & 0.0000 & 0.0078 & 0.5512 & 0.0000 & 0.0000 & - & 0.0075 & 0.0000 & 0.0117 & 0.0073 & 0.0038 & 0.0000 & 0.0000 & \textbf{0.0469}\\
\texttt{billy.talent} & \textbf{19.2817} & 0.0000 & 0.0000 & 0.0000 & 0.0000 & 0.0000 & 0.0988 & - & 0.0000 & 0.0000 & 0.0000 & \textbf{0.1107} & 0.0000 & 0.0000 & 0.0040\\
\texttt{the.killers} & 0.0000 & 0.0000 & 0.0000 & 0.0000 & 0.0000 & 0.0000 & 0.0277 & 0.0000 & - & 0.0000 & 0.0000 & 0.0277 & 0.0000 & 0.0000 & 0.0000\\
\texttt{the.beatles} & 0.0000 & 0.0000 & 0.0000 & 0.0000 & 0.0000 & 0.0000 & 0.0000 & 0.0000 & 0.0000 & - & 0.0000 & 0.0000 & 0.0000 & 0.0000 & 0.0000\\
\texttt{jack.johnson} & 0.0454 & 0.0000 & 0.0000 & 0.0038 & 0.0000 & 0.0039 & 0.0040 & 0.0000 & 0.0000 & 0.0000 & - & 0.0000 & 0.0000 & 0.0000 & 0.0000\\
\texttt{muse} & 0.0090 & 0.0090 & 0.0804 & 0.0000 & 0.0000 & 0.0474 & 0.0310 & 0.0226 & 0.0768 & 0.0000 & 0.0253 & - & \textbf{0.0483} & 0.0090 & 0.0345\\
\texttt{beatsteaks} & 0.0181 & 0.0151 &\textbf{0.1304} & 0.0030 & 0.0000 & 0.0362 & 0.0395 & 0.0339 &\textbf{0.2936} & 0.0000 & \textbf{0.2331} & 0.0362 & - & \textbf{0.1388} & 0.0000\\
\texttt{foo.fighters} & 0.0000 & \textbf{0.0392} & 0.0832 & 0.0030 & 0.0060 & 0.0084 & 0.0000 & 0.0169 & 0.0165 & 0.0000 & 0.0281 & 0.0090 & 0.0121 & - & 0.0029\\
\texttt{nirvana} & 0.0000 & 0.0000 & 0.0352 & 0.2014 & 0.0000 & 0.0000 & 0.0469 & 0.0113 & 0.0000 & \textbf{0.0156} & 0.0036 & 0.0038 & 0.0000 & 0.0000 & -\\\hline
SVM & \rotatebox[origin=c]{90}{\texttt{linkin.park}} & \rotatebox[origin=c]{90}{\texttt{coldplay}} & \rotatebox[origin=c]{90}{\mbox{ }\texttt{red.hot.chili.peppers}\mbox{\phantom{a} }} & \rotatebox[origin=c]{90}{\texttt{rammstein}} & \rotatebox[origin=c]{90}{\texttt{system.of.a.down}} & \rotatebox[origin=c]{90}{\texttt{metallica}} & \rotatebox[origin=c]{90}{\texttt{die.toten.hosen}} & \rotatebox[origin=c]{90}{\texttt{billy.talent}} & \rotatebox[origin=c]{90}{\texttt{the.killers}} & \rotatebox[origin=c]{90}{\texttt{the.beatles}} & \rotatebox[origin=c]{90}{\texttt{jack.johnson}} & \rotatebox[origin=c]{90}{\texttt{muse}} & \rotatebox[origin=c]{90}{\texttt{beatsteaks}} & \rotatebox[origin=c]{90}{\texttt{foo.fighters}} & \rotatebox[origin=c]{90}{\texttt{nirvana}}\\\hline
\texttt{linkin.park} & - & 0.0000 & 0.0000 & 0.0000 & 0.0000 & 0.0000 & 0.0000 & 0.0000 & 0.0000 & 0.0000 & 0.0000 & 0.0000 & 0.0000 & 0.0000 & 0.0000\\
\texttt{coldplay} & 0.0000 & - & 0.0000 & 0.0000 & 0.0000 & 0.0000 & 0.0000 & 0.0000 & 0.0000 & 0.0000 & 0.0000 & 0.0000 & 0.0000 & 0.0000 & 0.0000\\
\texttt{red.hot.chili.peppers} & 0.0000 & 0.0000 & - & 0.0000 & 0.0000 & 0.0000 & 0.0000 & 0.0000 & 0.0000 & 0.0000 & 0.0000 & 0.0000 & 0.0000 & 0.0000 & 0.0000\\
\texttt{rammstein} & 0.0000 & 0.0000 & 0.0000 & - & 0.0000 & 0.0000 & 0.0000 & 0.0000 & 0.0000 & 0.0000 & 0.0000 & 0.0000 & 0.0000 & 0.0000 & 0.0000\\
\texttt{system.of.a.down} & 0.0000 & 0.0000 & 0.0000 & 0.0000 & - & 0.0000 & 0.0000 & 0.0000 & 0.0000 & 0.0000 & 0.0000 & 0.0000 & 0.0000 & 0.0000 & 0.0000\\
\texttt{metallica} & 0.0000 & 0.0000 & 0.0000 & 0.0000 & 0.0000 & - & 0.0000 & 0.0000 & 0.0000 & 0.0000 & 0.0000 & 0.0000 & 0.0000 & 0.0000 & 0.0000\\
\texttt{die.toten.hosen} & 0.0000 & 0.0000 & 0.0000 & 0.0000 & 0.0000 & 0.0000 & - & 0.0000 & 0.0000 & 0.0000 & 0.0000 & 0.0000 & 0.0000 & 0.0000 & 0.0000\\
\texttt{billy.talent} & 0.0000 & 0.0000 & 0.0000 & 0.0000 & 0.0000 & 0.0000 & 0.0000 & - & 0.0000 & 0.0000 & 0.0000 & 0.0000 & 0.0000 & 0.0000 & 0.0000\\
\texttt{the.killers} & 0.0000 & 0.0000 & 0.0000 & 0.0000 & 0.0000 & 0.0000 & 0.0000 & 0.0000 & - & 0.0000 & 0.0000 & 0.0000 & 0.0000 & 0.0000 & 0.0000\\
\texttt{the.beatles} & 0.0000 & 0.0000 & 0.0000 & 0.0000 & 0.0000 & 0.0000 & 0.0000 & 0.0000 & 0.0000 & - & 0.0000 & 0.0000 & 0.0000 & 0.0000 & 0.0000\\
\texttt{jack.johnson} & 0.0000 & 0.0000 & 0.0000 & 0.0000 & 0.0000 & 0.0000 & 0.0000 & 0.0000 & 0.0000 & 0.0000 & - & 0.0000 & 0.0000 & 0.0000 & 0.0000\\
\texttt{muse} & 0.0000 & 0.0000 & 0.0000 & 0.0000 & 0.0000 & 0.0000 & 0.0000 & 0.0000 & 0.0000 & 0.0000 & 0.0000 & - & 0.0000 & 0.0000 & 0.0000\\
\texttt{beatsteaks} & 0.0000 & 0.0000 & 0.0000 & 0.0000 & 0.0000 & 0.0000 & 0.0000 & 0.0000 & 0.0000 & 0.0000 & 0.0000 & 0.0000 & - & 0.0000 & 0.0000\\
\texttt{foo.fighters} & 0.0000 & 0.0000 & 0.0000 & 0.0000 & 0.0000 & 0.0000 & 0.0000 & 0.0000 & 0.0000 & 0.0000 & 0.0000 & 0.0000 & 0.0000 & - & 0.0000\\
\texttt{nirvana} & 0.0000 & 0.0000 & 0.0000 & 0.0000 & 0.0000 & 0.0000 & 0.0000 & 0.0000 & 0.0000 & 0.0000 & 0.0000 & 0.0000 & 0.0000 & 0.0000 & - \\\hline
LR & \rotatebox[origin=c]{90}{\texttt{linkin.park}} & \rotatebox[origin=c]{90}{\texttt{coldplay}} & \rotatebox[origin=c]{90}{\mbox{ }\texttt{red.hot.chili.peppers}\mbox{\phantom{a}}} & \rotatebox[origin=c]{90}{\texttt{rammstein}} & \rotatebox[origin=c]{90}{\texttt{system.of.a.down}} & \rotatebox[origin=c]{90}{\texttt{metallica}} & \rotatebox[origin=c]{90}{\texttt{die.toten.hosen}} & \rotatebox[origin=c]{90}{\texttt{billy.talent}} & \rotatebox[origin=c]{90}{\texttt{the.killers}} & \rotatebox[origin=c]{90}{\texttt{the.beatles}} & \rotatebox[origin=c]{90}{\texttt{jack.johnson}} & \rotatebox[origin=c]{90}{\texttt{muse}} & \rotatebox[origin=c]{90}{\texttt{beatsteaks}} & \rotatebox[origin=c]{90}{\texttt{foo.fighters}} & \rotatebox[origin=c]{90}{\texttt{nirvana}}\\\hline
\texttt{linkin.park} & - & 0.0000 & 0.0000 & 0.0000 & 0.0000 & 0.0000 & 0.0000 & 0.0000 & 0.0000 & 0.0000 & 0.0000 & 0.0000 & 0.0000 & 0.0000 & 0.0000\\
\texttt{coldplay} & 0.0000 & - & 0.0000 & 0.0000 & 0.0000 & 0.0000 & 0.0000 & 0.0000 & 0.0000 & 0.0000 & 0.0000 & 0.0000 &\textbf{0.0030} & 0.0000 & 0.0000\\
\texttt{red.hot.chili.peppers} & 0.0000 & 0.0000 & - & 0.0000 & 0.0000 & 0.0000 & 0.0000 & 0.0000 & 0.0000 & 0.0078 & 0.0000 & 0.0000 & 0.0000 & 0.0000 & 0.0000\\
\texttt{rammstein} & 0.0000 & 0.0000 & 0.0000 & - & 0.0000 &\textbf{0.0056} & \textbf{0.0056} & 0.0000 & 0.0000 & 0.0000 & 0.0000 & 0.0000 & 0.0000 & 0.0000 & 0.0000\\
\texttt{system.of.a.down} & 0.0000 & 0.0000 & 0.0000 & 0.0965 & - & 0.0000 & 0.0000 & 0.0000 & 0.0000 & 0.0000 & 0.0000 & 0.0000 & 0.0000 & 0.0000 & 0.0000\\
\texttt{metallica} & 0.0000 & 0.0000 & 0.0000 & 0.0000 & 0.0000 & - & 0.0000 & 0.0000 & 0.0000 & \textbf{0.0117} & 0.0000 & 0.0000 & 0.0000 & 0.0000 & 0.0000\\
\texttt{die.toten.hosen} & 0.0000 & 0.0000 & 0.0000 & 0.0000 & 0.0000 & 0.0000 & - & 0.0000 & 0.0000 & 0.0000 & 0.0000 & 0.0000 & 0.0000 & 0.0000 & 0.0000\\
\texttt{billy.talent} & 0.0000 & 0.0000 & 0.0000 & 0.0000 & 0.0000 & 0.0000 & 0.0000 & - & 0.0000 & 0.0000 & 0.0000 & 0.0000 & 0.0000 & 0.0000 & 0.0000\\
\texttt{the.killers} & 0.0000 & 0.0000 & 0.0000 & 0.0000 & 0.0000 & 0.0000 & 0.0000 & 0.0000 & - & 0.0000 & 0.0000 & 0.0000 & 0.0000 & 0.0000 & 0.0000\\
\texttt{the.beatles} & 0.0000 & 0.0000 & 0.0000 & 0.0000 & 0.0000 & 0.0000 & 0.0000 & 0.0000 & 0.0000 & - & 0.0000 & 0.0000 & 0.0000 & 0.0000 & 0.0000\\
\texttt{jack.johnson} & 0.0000 & 0.0000 & 0.0000 & 0.0000 & 0.0000 & 0.0000 & 0.0000 & 0.0000 & 0.0000 & 0.0000 & - & 0.0000 & 0.0000 & 0.0000 & 0.0000\\
\texttt{muse} & 0.0000 & 0.0000 & 0.0000 & \textbf{0.1448} & 0.0000 & 0.0000 & 0.0000 & 0.0000 & 0.0000 & 0.0000 & 0.0000 & - & 0.0000 & 0.0000 & 0.0000\\
\texttt{beatsteaks} & 0.0000 & 0.0000 & 0.0000 & 0.0000 & 0.0000 & 0.0000 & 0.0000 & 0.0000 & \textbf{0.0027} & 0.0000 & 0.0000 & 0.0000 & - &\textbf{0.0030} & 0.0000\\
\texttt{foo.fighters} & 0.0000 & \textbf{0.0030} & 0.0000 & 0.0000 & 0.0000 & 0.0000 & 0.0000 & \textbf{0.0028} & 0.0000 & 0.0000 & 0.0000 &\textbf{0.0030} & \textbf{0.0030} & - & 0.0000\\
\texttt{nirvana} & 0.0000 & 0.0000 & 0.0000 & 0.0000 & 0.0000 & 0.0000 & 0.0000 & 0.0000 & 0.0000 & 0.0039 & 0.0000 & 0.0000 & 0.0000 & 0.0000 & - \\\hline
    \end{tabular}}
    \caption{Influence measure \eqref{eq:measure} using the constructed-by-design partition in the analysis of musical taste in Spotify.}
    \label{tab:musicdecades}
\end{table}

In view of Table~\ref{tab:musicdecades}, the most influential bands do not seem to coincide in the rankings resulting from using RF and LR as classifiers. However, we will use Pearson correlations on the numerical results as a measure of overall comparison of the rankings obtained (see Table \ref{cor:musicdecades}). The greatest similarities between rankings are found, in this order, when we study the influence of listening to The Killers, Coldplay, Foo Fighters, The Beatles and Die Toten Hosen. All of these bands have a correlation of around or above 0.6. 
\begin{table}[h!]
    \centering\resizebox{0.9\textwidth}{!}{
    \begin{tabular}{p{1.1 cm}p{1.1 cm}p{1.1 cm}p{1.1 cm}p{1.1 cm}p{1.3 cm}p{1.1 cm}p{1.1 cm}p{1.1 cm}p{1.1 cm}p{1.1 cm}p{1.3 cm}p{1.1 cm}p{1.1 cm}p{1.1 cm}}\hline
    \rotatebox[origin=c]{90}{\texttt{linkin.park}} & \rotatebox[origin=c]{90}{\texttt{coldplay}} & \rotatebox[origin=c]{90}{\mbox{ }\texttt{red.hot.chili.peppers}\mbox{ }} & \rotatebox[origin=c]{90}{\texttt{rammstein}} & \rotatebox[origin=c]{90}{\texttt{system.of.a.down}} & \rotatebox[origin=c]{90}{\texttt{metallica}} & \rotatebox[origin=c]{90}{\texttt{die.toten.hosen}} & \rotatebox[origin=c]{90}{\texttt{billy.talent}} & \rotatebox[origin=c]{90}{\texttt{the.killers}} & \rotatebox[origin=c]{90}{\texttt{the.beatles}} & \rotatebox[origin=c]{90}{\texttt{jack.johnson}} & \rotatebox[origin=c]{90}{\texttt{muse}} & \rotatebox[origin=c]{90}{\texttt{beatsteaks}} & \rotatebox[origin=c]{90}{\texttt{foo.fighters}} & \rotatebox[origin=c]{90}{\texttt{nirvana}}\\\hline
- &	0.8824 & - &	0.1451 &- &		-0.0099&	0.5984	&0.0448&	0.9290&	0.7143 &-&	-0.0822	&0.2272	&0.8704	&- \\\hline
     \end{tabular}}
    \caption{Correlations between the rankings obtained for RF and LR with CDP.}
    \label{cor:musicdecades}
\end{table}

\begin{table}[h!]
    \centering\resizebox{0.9\textwidth}{!}{
   \begin{tabular}{p{4 cm}|p{1.1 cm}p{1.1 cm}p{1.1 cm}p{1.1 cm}p{1.1 cm}p{1.1 cm}p{1.1 cm}p{1.1 cm}p{1.1 cm}p{1.1 cm}p{1.1 cm}p{1.1 cm}p{1.1 cm}p{1.1 cm}p{1.1 cm}}\hline
RF & \rotatebox[origin=c]{90}{\texttt{linkin.park}} & \rotatebox[origin=c]{90}{\texttt{coldplay}} & \rotatebox[origin=c]{90}{\mbox{ }\texttt{red.hot.chili.peppers}\mbox{\phantom{a}}} & \rotatebox[origin=c]{90}{\texttt{rammstein}} & \rotatebox[origin=c]{90}{\texttt{system.of.a.down}} & \rotatebox[origin=c]{90}{\texttt{metallica}} & \rotatebox[origin=c]{90}{\texttt{die.toten.hosen}} & \rotatebox[origin=c]{90}{\texttt{billy.talent}} & \rotatebox[origin=c]{90}{\texttt{the.killers}} & \rotatebox[origin=c]{90}{\texttt{the.beatles}} & \rotatebox[origin=c]{90}{\texttt{jack.johnson}} & \rotatebox[origin=c]{90}{\texttt{muse}} & \rotatebox[origin=c]{90}{\texttt{beatsteaks}} & \rotatebox[origin=c]{90}{\texttt{foo.fighters}} & \rotatebox[origin=c]{90}{\texttt{nirvana}}\\\hline
\texttt{linkin.park} &  - & 0.3989 & 0.0000 & 1.2427 & 1.1437 & 0.0031 & 0.0269 & 0.0062 & 0.0078 & 0.0000 & 0.0082 & 0.0138 & 0.0126 & 0.0000 & 0.0211\\
\texttt{coldplay} & 0.0000 &  - & 0.5779 & 0.0000 & 0.0000 & 0.0000 & 0.0000 & 0.0000 & 0.0000 & 0.0000 & 0.0000 & 0.0000 & 0.0000 & 0.0000 & 0.0000\\
\texttt{red.hot.chili.peppers} & 0.0093 &\textbf{0.8825} &  - & 0.3369 & 0.1669 & \textbf{0.4822} & 0.0000 & 0.0309 & 0.0117 & 0.0000 & 0.0408 & 0.0055 & 0.0067 & 0.0038 & 0.0070\\
\texttt{rammstein} & 0.0556 & 0.0328 & 0.0335 &  - &\textbf{1.4807} & 0.1669 & 0.0000 & 0.0340 & 0.0000 & 0.0000 & 0.0082 & 0.0083 & 0.0000 & 0.0000 & 0.0000\\
\texttt{system.of.a.down} & 0.0216 & 0.0383 & 0.1673 & \textbf{1.4683} &  - & 0.0433 & 0.0000 & 0.1453 & 0.0000 & 0.0000 & 0.0163 & 0.0083 & 0.0000 & 0.0000 & 0.0000\\
\texttt{metallica} & 0.0247 & 0.0301 & 0.0000 & 1.2674 & 1.1376 &  - & 0.0000 & 0.0309 & 0.0000 & 0.0000 & 0.0245 & 0.0083 & 0.0034 & 0.0038 & 0.0000\\
\texttt{die.toten.hosen} & 0.0309 & 0.0055 & 0.1673 & 1.2798 & 0.2844 & 0.0185 &  - & \textbf{0.1484} & 0.0156 & 0.0000 & 0.0109 & 0.0166 & 0.0000 & 0.0000 & \textbf{0.0352}\\
\texttt{billy.talent} & \textbf{21.3682} & 0.0820 & 0.0000 & 0.1329 & 0.2658 & 0.0464 & \textbf{0.0960} &  - & 0.0000 & 0.0000 & 0.0082 & \textbf{0.0775} & 0.0502 & 0.0000 & 0.0282\\
\texttt{the.killers} & 0.0000 & 0.0000 & 0.0000 & 0.0000 & 0.0000 & 0.0000 & 0.0307 & 0.0000 &  - & 0.0000 & 0.0000 & 0.0000 & 0.0000 & 0.0045 & 0.0000\\
\texttt{the.beatles} & 0.0000 & 0.0277 & 0.0608 & 0.0000 & 0.0000 & 0.0000 & 0.0000 & 0.0000 & 0.0000 &  - & 0.0000 & 0.0000 & 0.0000 & 0.0000 & 0.0000\\
\texttt{jack.johnson} & 0.0000 & 0.0000 & \textbf{1.3080} & 0.0037 & 0.0000 & 0.0000 & 0.0000 & 0.0000 & 0.0000 & 0.0000 &  - & 0.0000 & 0.0000 & \textbf{0.0115} & 0.0000\\
\texttt{muse} & 0.0042 & 0.0000 & 0.0078 & 0.0000 & 0.0000 & 0.0000 & 0.0192 & 0.0000 &\textbf{0.0165} & 0.0000 & 0.0000 &  - & 0.0000 & 0.0000 & 0.0000\\
\texttt{beatsteaks} & 0.0126 & 0.0000 & 0.0000 & 0.0000 & 0.0000 & 0.0000 & 0.0000 & 0.0045 & 0.0000 & 0.0000 &\textbf{0.3152} & 0.0000 &  - & 0.0000 & 0.0070\\
\texttt{foo.fighters} & 0.1004 & 0.0285 & 0.0044 & 0.0000 & 0.0000 & 0.0043 & 0.0000 & 0.0000 & 0.0041 & 0.0000 & 0.0162 & 0.0000 & \textbf{0.0753} &  - & 0.0070\\
\texttt{nirvana} & 0.1235 & 0.0000 & 0.0000 & 0.0000 & 0.0000 & 0.0000 & 0.0531 & 0.0000 & 0.0000 & 0.0000 & 0.0000 & 0.0000 & 0.0000 & 0.0000 & -\\\hline
SVM & \rotatebox[origin=c]{90}{\texttt{linkin.park}} & \rotatebox[origin=c]{90}{\texttt{coldplay}} & \rotatebox[origin=c]{90}{\mbox{ }\texttt{red.hot.chili.peppers}\mbox{\phantom{a}}} & \rotatebox[origin=c]{90}{\texttt{rammstein}} & \rotatebox[origin=c]{90}{\texttt{system.of.a.down}} & \rotatebox[origin=c]{90}{\texttt{metallica}} & \rotatebox[origin=c]{90}{\texttt{die.toten.hosen}} & \rotatebox[origin=c]{90}{\texttt{billy.talent}} & \rotatebox[origin=c]{90}{\texttt{the.killers}} & \rotatebox[origin=c]{90}{\texttt{the.beatles}} & \rotatebox[origin=c]{90}{\texttt{jack.johnson}} & \rotatebox[origin=c]{90}{\texttt{muse}} & \rotatebox[origin=c]{90}{\texttt{beatsteaks}} & \rotatebox[origin=c]{90}{\texttt{foo.fighters}} & \rotatebox[origin=c]{90}{\texttt{nirvana}}\\\hline
\texttt{linkin.park} &  - & 0.0000 & 0.0000 & 0.0000 & 0.0000 & 0.0000 & 0.0000 & 0.0000 & 0.0000 & 0.0000 & 0.0000 & 0.0000 & 0.0000 & 0.0000 & 0.0000\\
\texttt{coldplay} & 0.0000 &  - & 0.0000 & 0.0000 & 0.0000 & 0.0000 & 0.0000 & 0.0000 & 0.0000 & 0.0000 & 0.0000 & 0.0000 & 0.0000 & 0.0000 & 0.0000\\
\texttt{red.hot.chili.peppers} & 0.0000 & 0.0000 &  - & 0.0000 & 0.0000 & 0.0000 & 0.0000 & 0.0000 & 0.0000 & 0.0000 & 0.0000 & 0.0000 & 0.0000 & 0.0000 & 0.0000\\
\texttt{rammstein} & 0.0000 & 0.0000 & 0.0000 &  - & 0.0000 & 0.0000 & 0.0000 & 0.0000 & 0.0000 & 0.0000 & 0.0000 & 0.0000 & 0.0000 & 0.0000 & 0.0000\\
\texttt{system.of.a.down} & 0.0000 & 0.0000 & 0.0000 & 0.0000 &  - & 0.0000 & 0.0000 & 0.0000 & 0.0000 & 0.0000 & 0.0000 & 0.0000 & 0.0000 & 0.0000 & 0.0000\\
\texttt{metallica} & 0.0000 & 0.0000 & 0.0000 & 0.0000 & 0.0000 &  - & 0.0000 & 0.0000 & 0.0000 & 0.0000 & 0.0000 & 0.0000 & 0.0000 & 0.0000 & 0.0000\\
\texttt{die.toten.hosen} & 0.0000 & 0.0000 & 0.0000 & 0.0000 & 0.0000 & 0.0000 &  - & 0.0000 & 0.0000 & 0.0000 & 0.0000 & 0.0000 & 0.0000 & 0.0000 & 0.0000\\
\texttt{billy.talent} & 0.0000 & 0.0000 & 0.0000 & 0.0000 & 0.0000 & 0.0000 & 0.0000 &  - & 0.0000 & 0.0000 & 0.0000 & 0.0000 & 0.0000 & 0.0000 & 0.0000\\
\texttt{the.killers} & 0.0000 & 0.0000 & 0.0000 & 0.0000 & 0.0000 & 0.0000 & 0.0000 & 0.0000 &  - & 0.0000 & 0.0000 & 0.0000 & 0.0000 & 0.0000 & 0.0000\\
\texttt{the.beatles} & 0.0000 & 0.0000 & 0.0000 & 0.0000 & 0.0000 & 0.0000 & 0.0000 & 0.0000 & 0.0000 &  - & 0.0000 & 0.0000 & 0.0000 & 0.0000 & 0.0000\\
\texttt{jack.johnson} & 0.0000 & 0.0000 & 0.0000 & 0.0000 & 0.0000 & 0.0000 & 0.0000 & 0.0000 & 0.0000 & 0.0000 &  - & 0.0000 & 0.0000 & 0.0000 & 0.0000\\
\texttt{muse} & 0.0000 & 0.0000 & 0.0000 & 0.0000 & 0.0000 & 0.0000 & 0.0000 & 0.0000 & 0.0000 & 0.0000 & 0.0000 &  - & 0.0000 & 0.0000 & 0.0000\\
\texttt{beatsteaks} & 0.0000 & 0.0000 & 0.0000 & 0.0000 & 0.0000 & 0.0000 & 0.0000 & 0.0000 & 0.0000 & 0.0000 & 0.0000 & 0.0000 &  - & 0.0000 & 0.0000\\
\texttt{foo.fighters} & 0.0000 & 0.0000 & 0.0000 & 0.0000 & 0.0000 & 0.0000 & 0.0000 & 0.0000 & 0.0000 & 0.0000 & 0.0000 & 0.0000 & 0.0000 &  - & 0.0000\\
\texttt{nirvana} & 0.0000 & 0.0000 & 0.0000 & 0.0000 & 0.0000 & 0.0000 & 0.0000 & 0.0000 & 0.0000 & 0.0000 & 0.0000 & 0.0000 & 0.0000 & 0.0000 & -\\\hline
LR & \rotatebox[origin=c]{90}{\texttt{linkin.park}} & \rotatebox[origin=c]{90}{\texttt{coldplay}} & \rotatebox[origin=c]{90}{\mbox{ }\texttt{red.hot.chili.peppers}\mbox{\phantom{a}}} & \rotatebox[origin=c]{90}{\texttt{rammstein}} & \rotatebox[origin=c]{90}{\texttt{system.of.a.down}} & \rotatebox[origin=c]{90}{\texttt{metallica}} & \rotatebox[origin=c]{90}{\texttt{die.toten.hosen}} & \rotatebox[origin=c]{90}{\texttt{billy.talent}} & \rotatebox[origin=c]{90}{\texttt{the.killers}} & \rotatebox[origin=c]{90}{\texttt{the.beatles}} & \rotatebox[origin=c]{90}{\texttt{jack.johnson}} & \rotatebox[origin=c]{90}{\texttt{muse}} & \rotatebox[origin=c]{90}{\texttt{beatsteaks}} & \rotatebox[origin=c]{90}{\texttt{foo.fighters}} & \rotatebox[origin=c]{90}{\texttt{nirvana}}\\\hline
\texttt{linkin.park} &  - & 0.0000 & 0.0000 & 0.0185 & 0.0000 & 0.0000 & 0.0000 & 0.0000 & 0.0000 & 0.0000 & 0.0000 & 0.0000 & 0.0000 & 0.0000 & 0.0000\\
\texttt{coldplay} & 0.0000 &  - & 0.0532 & 0.0000 & 0.0000 & 0.0000 & 0.0000 & 0.0000 & 0.0000 & 0.0000 & 0.0000 & 0.0000 & 0.0000 & 0.0000 & 0.0000\\
\texttt{red.hot.chili.peppers} & 0.0000 & 0.0219 &  - & 0.0278 & 0.0031 & 0.0000 & 0.0000 & 0.0000 & 0.0000 & 0.0000 & 0.0000 & 0.0000 & 0.0000 & 0.0000 & 0.0000\\
\texttt{rammstein} &\textbf{0.0031} & 0.0000 & 0.0000 &  - & 0.0124 & \textbf{0.0093} & 0.0000 & 0.0000 & 0.0000 & 0.0000 & 0.0000 & 0.0000 & 0.0000 & 0.0000 & 0.0000\\
\texttt{system.of.a.down} & 0.0000 & 0.0000 & 0.0000 & \textbf{0.7913} &  - & 0.0000 & 0.0000 & 0.0000 & 0.0000 & 0.0000 & 0.0000 & 0.0000 & 0.0000 & 0.0000 & 0.0000\\
\texttt{metallica} & 0.0000 & 0.0027 & 0.0000 & 0.4884 & 0.0031 &  - & 0.0000 & 0.0000 & 0.0000 & 0.0000 & 0.0000 & 0.0000 & 0.0000 & 0.0000 & 0.0000\\
\texttt{die.toten.hosen} & 0.0000 & 0.0000 & 0.0000 & 0.3369 & 0.0031 & 0.0000 &  - & 0.0000 & 0.0000 & 0.0000 & 0.0000 & 0.0000 & 0.0000 & 0.0000 & 0.0000\\
\texttt{billy.talent} & 0.0000 & 0.0000 & 0.0000 & 0.0309 & \textbf{0.0402} & 0.0062 & 0.0000 &  - & 0.0000 & 0.0000 & 0.0000 & 0.0000 & 0.0000 & 0.0000 & 0.0000\\
\texttt{the.killers} & 0.0000 & 0.0000 & 0.0000 & 0.0000 & 0.0000 & 0.0000 & 0.0000 & 0.0000 &  - & 0.0000 & 0.0000 & 0.0000 & 0.0000 & 0.0000 & 0.0000\\
\texttt{the.beatles} & 0.0000 & 0.0000 & \textbf{0.0608} & 0.0000 & 0.0000 & 0.0000 & 0.0000 & 0.0000 & 0.0000 &  - & 0.0000 & 0.0000 & 0.0000 & 0.0000 & 0.0000\\
\texttt{jack.johnson} & 0.0000 & 0.0000 & 0.0456 & 0.0000 & 0.0000 & 0.0000 & \textbf{0.0036} & 0.0000 & 0.0000 & 0.0000 &  - & 0.0000 & 0.0000 & 0.0000 & 0.0000\\
\texttt{muse} & 0.0000 & 0.0000 & 0.0000 & 0.0000 & 0.0000 & 0.0000 & 0.0000 & 0.0000 & 0.0000 & 0.0000 & 0.0000 &  - & 0.0000 & 0.0000 & 0.0000\\
\texttt{beatsteaks} & 0.0000 & 0.0000 & 0.0000 & 0.0000 & 0.0000 & 0.0000 & 0.0000 & 0.0000 & 0.0000 & 0.0000 & 0.0000 & 0.0000 &  - & 0.0000 & 0.0000\\
\texttt{foo.fighters} & 0.0000 & \textbf{0.0285} & 0.0000 & 0.0000 & 0.0000 & 0.0000 & 0.0000 & 0.0000 & 0.0000 & 0.0000 & 0.0000 & 0.0000 & 0.0000 &  - & 0.0000\\
\texttt{nirvana} & 0.0000 & 0.0000 & 0.0000 & 0.0000 & 0.0000 & 0.0000 & 0.0000 & 0.0000 & 0.0000 & 0.0000 & 0.0000 & 0.0000 & \textbf{0.0067} & 0.0000 & -\\\hline
    \end{tabular}}
    \caption{Influence measure \eqref{eq:measure} using the hierarchical clustering partition in the analysis of musical taste in Spotify.}
    \label{tab:music_hcp}
\end{table}

Second, we do a similar study from the results obtained with the partitions prescribed by the hierarchical clustering. In this case, using RF as a classifier, Red Hot Chili Peppers is the most influential on listening to Coldplay and Metallica according to our influence measure; Rammstein, on System of a Down (and vice versa); Die Toten Hosen, to Billy Talent and Nirvana; Billy Talent, to Linking Park, Die Toten Hosen and Muse; Jack Johnson, to Red Hot Chilly Peppers and Foo Fighters; Muse, to The Killers; Beatsteaks, to Jack Johnson; and Foo Fighters, to Beatsteaks. Regarding the results based on LR, we mention the following issues. Our influence measure shows that Rammstein is, in this case, the most influential band in listening to Linkin Park and Metallica; System of a Down, to Rammsten; Billy Talent, to Systems of a Down; The Beatles, to Red Hot Chilly Peppers; Jack Johnson, to Die Toten Hosen; Foo Fighters, to Coldplay, and finally, Nirvana, to Beatsteaks.

Table \ref{cor:music_hcp} shows the Pearson correlations between the influence measures obtained for each band when using RF and LR as classifiers. In this case, the greatest similarities between the rankings are found, in this order, when we examine the influence of listening to Rammstein, Red Hot Chilli Peppers, and Coldplay. 

\begin{table}[h!]
    \centering\resizebox{0.9\textwidth}{!}{
    \begin{tabular}{p{1.3 cm}p{1.1 cm}p{1.1 cm}p{1.1 cm}p{1.1 cm}p{1.1 cm}p{1.3 cm}p{1.1 cm}p{1.1 cm}p{1.1 cm}p{1.1 cm}p{1.3 cm}p{1.3 cm}p{1.1 cm}p{1.1 cm}}\hline
     \rotatebox[origin=c]{90}{\texttt{linkin.park}} & \rotatebox[origin=c]{90}{\texttt{coldplay}} & \rotatebox[origin=c]{90}{\mbox{ }\texttt{red.hot.chili.peppers}\mbox{ }} & \rotatebox[origin=c]{90}{\texttt{rammstein}} & \rotatebox[origin=c]{90}{\texttt{system.of.a.down}} & \rotatebox[origin=c]{90}{\texttt{metallica}} & \rotatebox[origin=c]{90}{\texttt{die.toten.hosen}} & \rotatebox[origin=c]{90}{\texttt{billy.talent}} & \rotatebox[origin=c]{90}{\texttt{the.killers}} & \rotatebox[origin=c]{90}{\texttt{the.beatles}} & \rotatebox[origin=c]{90}{\texttt{jack.johnson}} & \rotatebox[origin=c]{90}{\texttt{muse}} & \rotatebox[origin=c]{90}{\texttt{beatsteaks}} & \rotatebox[origin=c]{90}{\texttt{foo.fighters}} & \rotatebox[origin=c]{90}{\texttt{nirvana}}\\\hline
-0.0756&	0.4911&	0.6255&	0.8122&	0.1923&	0.2026&	-0.1640&- 	&	- 	&	- 	&	- 	&	- 	&	-0.1328 &	- 	&	-  	\\\hline
     \end{tabular}}
    \caption{Correlations between the rankings obtained for RF and LR with HCP.}
    \label{cor:music_hcp}
\end{table}

Finally, one aspect to be considered is the possible effect of partitioning on the resulting ranking. Although some similarities can apparently be observed between the rankings obtained under CDP and HCP, we determine the corresponding Pearson correlations in Table \ref{cor:comparison} between the rankings obtained when using RF as classifier as well as the case of considering LR.

\begin{table}[h!]
    \centering\resizebox{0.9\textwidth}{!}{
    \begin{tabular}{p{0.6 cm}|p{1.2 cm}p{1.2 cm}p{1.2 cm}p{1.2 cm}p{1.2 cm}p{1.2 cm}p{1.3 cm}p{1.2 cm}p{1.2 cm}p{1.2 cm}p{1.2cm}p{1.3 cm}p{1.3 cm}p{1.2 cm}p{1.2 cm}}\hline
    &  \rotatebox[origin=c]{90}{\texttt{linkin.park}} & \rotatebox[origin=c]{90}{\texttt{coldplay}} & \rotatebox[origin=c]{90}{\mbox{ }\texttt{red.hot.chili.peppers}\mbox{ }} & \rotatebox[origin=c]{90}{\texttt{rammstein}} & \rotatebox[origin=c]{90}{\texttt{system.of.a.down}} & \rotatebox[origin=c]{90}{\texttt{metallica}} & \rotatebox[origin=c]{90}{\texttt{die.toten.hosen}} & \rotatebox[origin=c]{90}{\texttt{billy.talent}} & \rotatebox[origin=c]{90}{\texttt{the.killers}} & \rotatebox[origin=c]{90}{\texttt{the.beatles}} & \rotatebox[origin=c]{90}{\texttt{jack.johnson}} & \rotatebox[origin=c]{90}{\texttt{muse}} & \rotatebox[origin=c]{90}{\texttt{beatsteaks}} & \rotatebox[origin=c]{90}{\texttt{foo.fighters}} & \rotatebox[origin=c]{90}{\texttt{nirvana}}\\\hline
RF & 0.9999 & 	-0.1964&	-0.1546&	{0.9484}&	{0.7776}&	-0.1768&	0.0469&	-0.2691&	-0.1001& -&	0.9859&	{0.8619}&	-0.0137&	-0.2259& 0.3248 \\
LR & - &	{0.7746}& -	&0.3316&-	&	{0.8190}&	-0.0769&-&-&-&-&-& -0.1132&-&-\\\hline
     \end{tabular}}
    \caption{Correlations between the rankings obtained under CDP and HCP.}
    \label{cor:comparison}
\end{table}

When using RF, the highest similarities in the rankings are found for Linkin Park, Jack Johnson, Rammstein, Muse and System of a Down, with correlations above 0.75. In these cases, the use of different partitions on the features causes small numerical variations in the rankings. This is not common for the rest of the bands, where even negative correlations are obtained. Nor is it when using LR as a classifier, where only Metallica and Coldplay exceed the 0.75 correlation threshold.

\section{Concluding remarks}\label{sec:conclusions}

In this paper, we have introduced a method to analyze the influence of specific features in a classification problem where dependencies might exist between them. Our influence measure extends \citeauthor{Datta2015}'s approach to accommodate the existence of a coalitional structure on such attributes. Through properties typical in cooperative game theory, we provide an axiomatic characterization of this measure. Furthermore, we demonstrate that the Banzhaf-Owen value serves as a special case of our measure under conditions of binary databases and no repeated profiles. The methodology presented here is subsequently applied to investigate influential variables in traffic accidents or to analyze musical taste on Spotify.

In this context, several lines of research stemming from our proposal remain open and need to be addressed in the future. From a practical perspective, the methodology presented in this paper can be used in a wide range of real-world applications. One example is the booming sector of streaming platforms such as Amazon Prime, Disney+, Max, or Netflix. It would be particularly interesting to conduct a similar analysis to what we performed for the Spotify database, aiming to identify which audiovisual content has the greatest (or least) impact on the viewing of others. The inherent similarities among some of these media resources make our influence measure with dependencies applicable in this domain. This would enable us to propose personalized recommendations to users, which is of major relevance to platforms' stakeholders \citep{bourreau2022}. However, significant computational challenges arise when applying our methodology to real-world problems due to the sheer and diverse amount of information available in today's society. 
Besides, we assumed in Section~\ref{sec:preli} that each element of a sample $\mathcal{M}$ appears exactly once. Alternatively, we could consider that each of the $n$ distinct profiles in the sample, $(X^i,Y^i)$, with $i\in N$, has an absolute frequency of $n(i)$, resulting in a total sample size of $\sum_{i\in N}n(i)$. By extending the properties used in Theorem~\ref{th:1}, particularly those of disjoint unions and symmetry, to this context, we can characterize a weighted influence measure that accommodates dependent features and accounts for repeated profiles. Specifically, for each $l\in K$, this measure is defined by 
\begin{equation*}
    \hat{I}_l(X, Y, P, \mathcal{M})= C \cdot \hat{\Psi}_l(X, Y, P, \mathcal{M}) = C \cdot 
\sum_{\substack{i\in N\,:
\\ (X^i, Y^i)\in \mathcal{M}^t
}
}
\, 
{\sum_{\substack{((X^i_{-l},\,a_{l}),\, b)\in \mathcal{M}^t\,:
\\ a_{l}\in  \mathcal{A}_{l},\,\, b\in \mathcal{B}
}
}
}n(i)\, \cdot\, |Y^i-b|.
\end{equation*}
Even so, the exact calculation in practice of our proposed influence measure reaches exponential computational complexity as does the Banzhaf-Owen value for a sufficiently large number of players involved (see \citealp{deng1994complexity}). To overcome this issue, it seems intuitive that our Banzhaf-Owen value-based influence measure can also be approximated by statistical sampling in large-scale settings, following the ideas in \cite{Saavedra-nieves2021} for estimating such a game-theoretical solution concept. Finally, we would like to highlight the potential value of conducting a more comprehensive study on the performance of our influence measure using other classifiers in the existing literature. To the best of the authors' knowledge, random forests, support vector machines, and logistic regression-type classifiers are among the most widely used methods for binary data. We have, however, verified the strong dependence of the chosen databases and the response variable determined by a certain classifier. Furthermore, we took base implementations of these methods in the \texttt{RWeka} library. Calibrating different parameters of these estimators, as well as testing other classification methods, might yield more conclusive results for some of the analyzed scenarios.
%\textcolor{blue}{Finally, we would like to note that a comprehensive study of the performance in practice of the wide variety of classifiers existing in the literature would be of interest. To the best of the authors' knowledge, we have focused on only three of the most widely used classifiers and we have verified the strong dependence of the database on the classifier used for the determination of the response variable.} 

\section*{Acknowledgments}
This work is part of the R+D+I project PID2021-124030NB-C32, granted by MICIU/AEI/

\noindent 10.13039/501100011033/ and by ``ERDF A way of making Europe''/EU. This research was also funded by {\it Grupos de Referencia Competitiva} ED431C 2021/24 from the \emph{Conseller{\'i}a de Cultura, Educaci{\'o}n e Universidades, Xunta de Galicia}.
Authors also thank the computational resources of the {\it Centro de Supercomputaci{\'o}n de Galicia} (CESGA).

\subsection*{Declaration of interest}
The authors declare that there is no conflict of interest.

\renewcommand{\BBAA}{and}
\bibliographystyle{apacite}
\bibliography{referenciassinDOI}

%\clearpage
%\appendix
%
%\setcounter{page}{1}
%\setcounter{table}{0}
%\renewcommand{\thetable}{\Alph{section}.\arabic{table}}
%\setcounter{figure}{0}
%\renewcommand\thefigure{\thesection.\arabic{figure}}  
%{\color{black}

\end{document}